\documentclass[12pt]{article}
\usepackage{setspace}
\usepackage[margin=1.25in]{geometry}
 
\usepackage{amsmath,amsthm,amssymb}
\usepackage{verbatim}
\usepackage{relsize}
\usepackage{listings}
\usepackage{xcolor}
\usepackage{hyperref}
\usepackage{graphicx}
\graphicspath{{./Images/}}

\usepackage{tikz}
\usetikzlibrary{arrows}

\newcommand{\PP}{\mathbb{P}}
\DeclareMathOperator{\E}{\mathbb{E}}
\DeclareMathOperator{\R}{\mathbb{R}}

\newcommand{\indep}{\perp \!\!\! \perp}

\usepackage{amsmath}

\DeclareMathOperator*{\argmin}{arg\,min}

\usepackage{mathtools}

\usepackage{tikz}
\usetikzlibrary{arrows}

\newcommand{\appendixhead}%
{\textbf{\huge Appendices}
\vspace{0.25in}}

\newtheorem{assumption}{Assumption}
\newtheorem{theorem}{Theorem}

\newtheorem{corollary}{Corollary}[theorem]
\newtheorem{remark}{Remark}

\usepackage{subcaption}

\title{Bias Formulas for Violations of Proximal Identification Assumptions}
\author{Raluca Cobzaru, Roy Welsch, Stan Finkelstein, Kenney Ng, Zach Shahn}

\begin{document}

\maketitle

\begin{singlespace}
\begin{abstract}
Causal inference from observational data often rests on the unverifiable assumption of no unmeasured confounding. Recently, Tchetgen Tchetgen and colleagues have introduced proximal inference to leverage negative control outcomes and exposures as proxies to adjust for bias from unmeasured confounding \cite{tchetgen_introduction_2020}. However, some of the key assumptions that proximal inference relies on are themselves empirically untestable. Additionally, the impact of violations of proximal inference assumptions on the bias of effect estimates is not well understood. In this paper, we derive bias formulas for proximal inference estimators under a linear structural equation model data generating process. These results are a first step toward sensitivity analysis and quantitative bias analysis of proximal inference estimators. While limited to a particular family of data generating processes, our results may offer some more general insight into the behavior of proximal inference estimators.
\end{abstract}
\end{singlespace}

\section{Introduction}

Causal inference using observational data often rests on the assumption of no unmeasured confounding. This assumption is not empirically verifiable, but sensitivity analysis methods (e.g., \cite{cornfield_smoking_2009, ding_sensitivity_2016, robins_sensitivity_2000, rosenbaum_assessing_1983}) are available to assess robustness of results to possible unmeasured confounding. Alternatively, investigators might turn to methods such as instrumental variable analysis or difference-in-differences, which depend on different assumptions. Sensitivity analyses for violations of the assumptions required by these alternative methods are also available (\cite{brookhart_preference-based_2007, rambachan_more_2022}). 

There has been recent interest in the use of negative control methods to detect and resolve confounder bias. A negative control outcome (NCO) is a variable known not to be causally affected by the treatment of interest, while a negative control exposure (NCE) is a variable known not to causally affect the outcome of interest \cite{shi_selective_2020}. Tchetgen Tchetgen and colleagues have developed a Proximal Inference framework (\cite{cui_semiparametric_2020, miao_identifying_2018, tchetgen_introduction_2020}) which uses NCE-NCO pairs sharing the same unmeasured confounders as the treatment-outcome relationship of interest as proxies to adjust for unmeasured confounding. 

However, some of the assumptions that proximal inference relies on are themselves empirically untestable \cite{tchetgen_introduction_2020} and bias resulting from violations of proximal identification assumptions is not fully understood. In this paper, we characterize bias from violations of proximal inference assumptions in a linear structural equation data generating process. Our results build understanding of the sensitivity of proximal inference to assumption violations and serve as a first step toward sensitivity analysis and quantitative bias analysis \cite{lash_good_2014} tools for proximal inference.

The organization of the paper is as follows. In Section 2, we review proximal inference. In Section 3, we describe the forms of bias that we will study. In Sections 4 and 5, we derive bias formulas for the settings described in Section 3. In Section 6, we present numerical experiments based on the bias formulas from Sections 4 and 5 to explore their implications. In Section 7, we conclude by discussing some potential insights into the sensitivity of proximal inference estimators gained from our results. \color{black}

\section{Proximal Identification of the Average Treatment Effect}

\subsection{Review of Definitions and Assumptions}

We use the potential outcome framework \cite{rubin_estimating_1974} to formally define causal effects. Let $A$ denote the binary treatment of interest, $Y$ the observed post-treatment outcome, and $Y(a), a = 0, 1$ the potential (counterfactual) outcome that would have been observed had treatment $A$ been set to $a$. We implicitly make the no-interference assumption that the potential outcome of each individual does not depend on the treatments received by other individuals \cite{vanderweele_bias_2011}. We aim to estimate the average causal effect (ACE) of $A$ on $Y$, defined as $\psi = \E[Y(1) - Y(0)]$.

Let $L$ denote the set of measured covariates. We make the standard assumptions of Consistency and Positivity, defined below.
\begin{assumption}{\textnormal{(Consistency)}} \label{assumption: consistency}
$Y = Y(A)$ 
\end{assumption}
In other words, the observed value of $Y$ under treatment $A$ coincides with the counterfactual outcome that would have been observed under the same treatment value. Thus, we only observe the counterfactual outcome corresponding to the treatment value that was actually administered in our data.

\begin{assumption}{\textnormal{(Positivity)}} \label{assumption: positivity}
$0 < \PP(A = a | L) < 1$, for $a = 0, 1$
\end{assumption}
Assumption 2 states that both exposure levels are observed at all levels of the observed covariates $L$. 

Many analyses further make the assumption that there is no unobserved confounding, i.e. that observed covariates block all non-direct (or `backdoor') causal paths between treatment and outcome. 
\begin{assumption}{\textnormal{(Exchangeability)}} \label{assumption: exchangeability}
$Y(a) \indep A \ | \ L$, for $a = 0, 1$
\end{assumption}

Under Assumptions \ref{assumption: consistency}-\ref{assumption: exchangeability}, the counterfactual mean $\E[Y(a)]$ is identified by the g-formula (introduced in \cite{greenland_identifiability_1986}):
\begin{equation} \label{eq: g-formula}
    \E[Y(a)] = \sum_{l} \E[Y| A = a, L = l] \PP(L = l).
\end{equation}

Exchangeability is a strong assumption that is empirically untestable. \cite{miao_confounding_2020} propose an alternative to Assumption \ref{assumption: exchangeability} that allows us to identify the counterfactual mean $\E[Y(a)]$ despite the presence of unobserved confounding. We review the alternative conditions developed by \cite{miao_confounding_2020} leading to the \textit{proximal g-formula}, a counterpart to \eqref{eq: g-formula} allowing for some unobserved confounding. 

\begin{figure}[!h]
    \centering
    \includegraphics[width=0.27\linewidth]{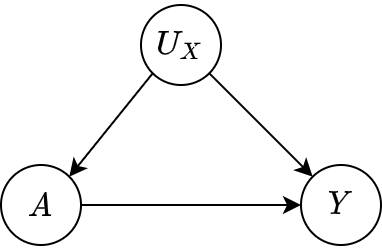}\quad\quad
    \includegraphics[width=0.27\linewidth]{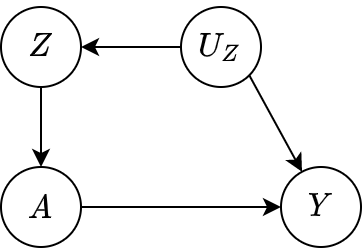}\quad\quad
    \includegraphics[width=0.27\linewidth]{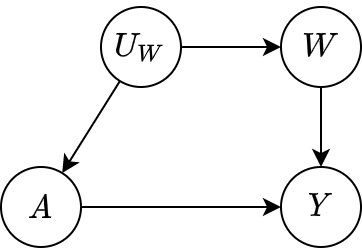}
    \caption{DAGs representing the three types of variables $(X,Z,W)$ partitioning $L$}
    \label{fig: L-types}
\end{figure}

As in \cite{cui_semiparametric_2020}, we consider a (potentially multidimensional) variable $L$ that can be partitioned into three types of variables $(X,Z,W)$, such that 
\begin{itemize}
    \item [1)] $X$ includes observed variables that may be common causes of $A$ and $Y$ (observed confounders)
    \item [2)] $Z$ includes treatment-inducing confounding proxies, i.e. $Z$ includes causes of $A$ that share an unmeasured common cause $U_Z$ with $Y$
    \item [3)] $W$ includes outcome-inducing confounding proxies, i.e. $W$ includes causes of $Y$ that share an unmeasured common cause $U_W$ with $A$
\end{itemize}
Figure \ref{fig: L-types} contains DAGs representing each of the proxy types included in $L$. 

In \cite{miao_confounding_2020}, exchangeability is replaced with the assumptions: 
\begin{assumption}{\textnormal{(Treatment-inducing confounding proxy)}}\label{assumption: NCE}
\begin{equation}
    Y(a,z) = Y(a), \ \ \textnormal{for all } a, z \label{eq: NCE-assumption}
\end{equation}
\end{assumption}
\begin{assumption}{\textnormal{(Outcome-inducing confounding proxy)}}\label{assumption: NCO}
\begin{equation}
    W(a,z) = W, \ \ \textnormal{for all } a, z \label{eq: NCO-assumption}
\end{equation}
\end{assumption}
\begin{assumption}{\textnormal{(Latent unconfoundedness)}}\label{assumption: latent} If $U$ denotes the set of unobserved confounders, then
\begin{align}
    Z &\indep (Y(a),W) \ | \ U,X \label{eq: latent-nce}\\
    W &\indep A \ | \ U, X \label{eq: latent-nco}
\end{align}
\end{assumption}
Assumption \ref{assumption: NCE} states that $Z$ does not have a direct effect on $Y$ upon intervening on $A$, while Assumption \ref{assumption: NCO} states that neither $A$ nor $Z$ have a causal effect on $W$. Past works \cite{shi_selective_2020} refer to variables $Z$ satisfying \eqref{eq: NCE-assumption} and \eqref{eq: latent-nce} as negative control exposure (NCE) variables, and to variables $W$ satisfying \eqref{eq: NCO-assumption} and \eqref{eq: latent-nco} as negative control outcome (NCO) variables. This terminology is based on negative control methods employing variables that share a confounding mechanism with the treatment-outcome relationship in view to detect bias in epidemiological research. Although there is a subtle distinction between the proxy and negative control nomenclature when discussing the design of observational studies \cite{tchetgen_introduction_2020}, for the theoretical analysis employed in this paper we will be using \textit{treatment-inducing} (\textit{outcome-inducing}) \textit{confounding proxies} and \textit{NCE} (\textit{NCO}) variables interchangeably. 

In addition to Assumptions \ref{assumption: consistency}-\ref{assumption: latent}, \cite{miao_identifying_2018} introduce the following \textit{completeness} conditions for the identification of $\E[Y(a)]$:
\begin{assumption}{\textnormal{(Completeness)}}\label{assumption: Completeness}
For any $a$, $x$ and for any square-integrable function $g$:
\begin{itemize}
    \item [(a)] If $\E[g(U) | Z, A = a, X = x] = 0$ almost surely, then $g(U) = 0$ almost surely.
    \item [(b)] If $\E[g(Z) | W, A = a, X = x] = 0$ almost surely, then $g(Z) = 0$ almost surely.
\end{itemize}
\end{assumption}

Assumption \ref{assumption: Completeness}(a) can be interpreted as a requirement that the NCE $Z$ has enough variability relative to the variability of $U$; similarly, assumption \ref{assumption: Completeness}(b) requires the variability of $W$ to be large enough relative to the variability of $Z$. Under conditions \ref{assumption: Completeness}(a) and (b), we can essentially account for $U$ in our ACE estimate without either measuring or modeling the distribution of $U$. The role of completeness will be further explored in Section \ref{sec: proximal estimation}, where we outline the analytical framework by which the ACE is estimated using the proximal g-formula. 

Completeness assumption \ref{assumption: Completeness}(a) has a simple interpretation in the case where confounders $U$ and the negative control pair $(Z,W)$ are all categorical. As mentioned in \cite{cui_semiparametric_2020}, if $(U,Z,W)$ are categorical with respective number of categories $(d_u, d_z, d_w)$, then completeness \ref{assumption: Completeness}(a) requires that
\begin{equation}
    \min (d_z, d_w) \ge d_u    
\end{equation}
In other words, proximal inference can account for unmeasured confounding if the number of categories of $U$ is less than that of either $Z$ or $W$. This leads to the practical recommendation to measure a rich set of baseline characteristics (which can be used as negative controls), such that the proximal identification approach has a higher chance of mitigating unmeasured confounder bias \cite{cui_semiparametric_2020}. There is not such a straightforward method for expressing the completeness condition in the case of continuous $U$ and negative controls $(Z,W)$, though some intuition for nonparametric regression results from \cite{cui_semiparametric_2020}. In Section \ref{sec: bias-fmlas}, we investigate the behavior of proximal inference in LSEM setups in which the completeness assumption \ref{assumption: Completeness}(a) is violated.

Lastly, to be valid proxies the variables $(Z,W)$ must be $U-relevant$:
\begin{assumption}{\textnormal{($U$-relevance)}}\label{assumption: relevance}
\begin{align}
    Z &\not \indep U \ | \ A, X \\
    W &\not \indep U \ | \ X
\end{align}
\end{assumption}
The $U$-relevance assumption (also known as $U$-comparability \cite{shi_selective_2020}) requires the unmeasured confounders $U$ of the $A-Y$ relationship to be the same as the unmeasured confounders of the $A-W$ and $Z-Y$ secondary treatment-outcome associations. This is such that, by the negative control framework, any non-null $A-W$ or $Z-Y$ association can be attributed to $U$ confounding the $A-Y$ relationship (while null associations imply no empirical evidence of unmeasured confounding).

Throughout this paper, we suppress the observed confounders $X$ unless otherwise stated. While we do not include $X$ in the sensitivity analysis discussion of Section \ref{sec: numerical-experiments}, the addition of $X$ is a straightforward extension of our bias formulations.


\subsection{Estimating the Proximal g-Formula via Moment Restriction} \label{sec: proximal estimation}

\cite{miao_confounding_2020} introduce the notion of an \textit{outcome confounding bridge function}, which transforms the negative control outcome $W$ to match the confounding effect of $U$ on $Y$. More precisely, an outcome confounding bridge function $h(W, A, X)$ is a function satisfying
\begin{equation}\label{eq: outcome-bridge-def}
    \E[Y | U, A = a, X = x] = \E[h(W, A, X) | U, A = a, X = x]
\end{equation}
for all values of $a, x$. In other words, if function $h(W, A, X)$ exists, then the confounding effect of $U$ on the transformed variable $h(W, a, X)$ equals the confounding effect of $U$ on $Y$ at exposure level $A = a$. Given assumptions \ref{assumption: consistency}, \ref{assumption: NCO}, \ref{assumption: latent}, and \ref{assumption: relevance}, \cite{miao_confounding_2020} infer that
\begin{equation}\label{eq: outcome-bridge-ace-overall}
    \E[Y(a)] = \E[h(W, a, X)] \textnormal{ for all $a = 0, 1$}
\end{equation}
which means $\E[Y(a)]$ can be estimated following the identification of an outcome bridge function $h(W, A, X)$, if such a function is assumed to exist.

\cite{cui_semiparametric_2020, miao_identifying_2018} established the following proximal identification result for the outcome confounding bridge function that leverages the distribution of a NCE $Z$:

\begin{theorem}\label{theorem: outcome-bridge}
Suppose there exists an outcome confounding bridge function $h(w,a,x)$ solving the Fredholm integral equation 
\begin{equation} \label{eq: outcome-bridge-Z}
    \E[Y| Z, A, X] = \int h(w, A, X) dF(w | Z, A, X)
\end{equation}
almost surely. Then, under Assumptions \ref{assumption: consistency}, \ref{assumption: positivity}, \ref{assumption: NCE}-\ref{assumption: latent}, and \ref{assumption: Completeness}(a),
\begin{equation} \label{eq: outcome-bridge-U}
    \E[Y | U, A, X] = \int h(w, A, X) dF(w | U, X)
\end{equation}
almost surely.
\end{theorem}

Under Assumption \ref{assumption: latent}, we have $\E[Y(a)] = \E\left[\E\left[Y| U, A = a, X = x \right] \right]$ for all $a, x$. The counterfactual mean $\E[Y(a)]$ can then be computed as follows:
\begin{corollary}{\textnormal{(Proximal g-formula)}}
If \eqref{eq: outcome-bridge-U} holds almost surely, then the counterfactual mean $\E[Y(a)]$, $a = 0, 1$ is nonparametrically identified by
\begin{equation} \label{eq: proximal-g-fmla}
    \E[Y(a)] = \int_\mathcal{X} \int h(w, a, x) dF(w | x) dF(x)
\end{equation}
and the ACE is identified by
\begin{equation} \label{eq: proximal ace}
    \psi = \int_\mathcal{X} \int \left\{h(w, 1, x) - h(w, 0, x) \right\} dF(w | x) dF(x) 
\end{equation}
\end{corollary}

\begin{remark}\label{remark1} \cite{cui_semiparametric_2020} establish a similar proximal identification result for the existence and identification of a treatment confounding bridge function $q(Z, A, X)$ that leverages the NCO variable $W$ (and an assumption analogous to completeness Assumption \ref{assumption: Completeness}) instead. Due to the higher complexity of $\frac{1}{\PP(A = a | U, X)}$ relative to $\E[Y | U, A, X]$ in our chosen LSEMs, we delegate sensitivity analysis involving the treatment confounding bridge function to future work. 
\end{remark}

Assuming the outcome confounding bridge function $h(W, A, X)$ exists and is identifiable as a solution to \eqref{eq: outcome-bridge-U}, \cite{tchetgen_introduction_2020, miao_confounding_2020} provide a practical approach for estimating the proximal g-formula using the generalized method of moments (GMM). Suppose one has access to $n$ i.i.d. samples $D_i = (A_i, Y_i, L_i)$, $L_i = (X_i, Z_i, W_i)$ (where $Z$, $W$ are assumed to be correctly classified as treatment- and outcome-inducing confounding proxies, respectively). Moreover, suppose one has specified a parametric model for the confounding bridge, $h(W, A, X) = h(W, A, X; b)$ (e.g., $h(W, A, X; b)$ is linear in $W$, $A$, $X$ with unknown parameter $b$). The true model for $h(W, A, X)$ is unknown, but one may fit a fairly flexible model (including, for instance, splines or interaction terms) to obtain a reasonable estimation in practice. 

We define the target parameter $\theta = (b, \psi)$ to encode the parameters $b$ of $h(W,A,X;b)$ and the ACE $\psi$, along with the moment restrictions
\begin{equation}
    h(D_i; \theta) = \begin{pmatrix} \left\{Y_i - h(W_i, A_i, X_i) \right\} \times \begin{pmatrix} 1 & A_i & Z_i & X_i & A_i X_i & A_i Z_i \end{pmatrix}^T \\ 
    \psi - \left\{h(W_i, V_i, 1; b) - h(W_i, V_i, 0; b) \right\} \end{pmatrix}
\end{equation}
Then, if $m_n(\theta) = \frac{1}{n} \sum_{i=1}^n h(D_i; \theta)$, the GMM solves
\begin{equation}
    \hat \theta = \argmin_\theta m_n^T(\theta) m_n(\theta)
\end{equation}
which can be equivalently written as
\begin{equation}
    \hat \theta = \begin{pmatrix} \hat b \\ \hat \psi \end{pmatrix} = \begin{pmatrix} \argmin_b \E\left[Y - h(W_i, A_i, X_i) | Z_i, A_i, X_i \right]^2 \\ 
    \argmin_{\psi} \left(\psi - \left\{h(W_i, V_i, 1; b) - h(W_i, V_i, 0; b) \right\} \right)^2 \end{pmatrix}
\end{equation}
The resulting parameter estimate $\hat b$ is unbiased by \eqref{eq: outcome-bridge-def}, while the ACE estimate $\hat \psi$ is unbiased by \eqref{eq: outcome-bridge-ace-overall}. 

\section{Bias Settings}

We have so far collected a series of untestable assumptions 4-8 that replace exchangeability and account for the effect of unmeasured confounders $U$ without directly modeling or estimating $U$. 
The impact on the direction and/or magnitude of bias resulting from violations of these assumptions has not been explored. We trust the analyst to identify true negative control exposures and outcomes in this work (Assumptions 4 and 5), as subject matter knowledge should often be quite reliable on this point. Latent unconfoundedness (Assumption 6) is not really an assumption since it presumably holds for some sufficiently rich $U$. But the richer (or higher dimensional) the $U$ required to satisfy latent unconfoundedness Assumption 6, the less plausible it is that completeness (Assumption 7) or $U$-relevance (Assumption 8) hold. If many components of $U$ are common causes of the negative control exposures and outcomes, then completeness (Assumption 7) is difficult to satisfy. And if many components of $U$ are required to block all backdoor paths between A and Y, then they are less likely to all be associated with both $Z$ and $W$, violating Assumption 8. 

\begin{figure}[!h]
    \centering
    \includegraphics[width=0.57\linewidth]{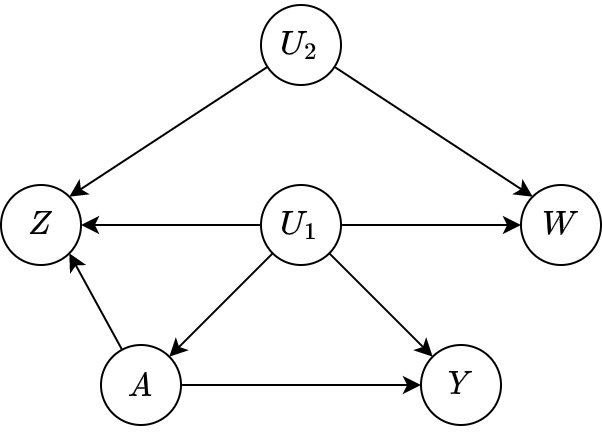}
    \caption{DAG encoding the causal relationships among variables in \eqref{eq: U2-ZW} in which completeness 7(a) is violated}
    \label{fig: completeness}
\end{figure}

\begin{figure}[!h]
    \centering
    \includegraphics[width=0.57\linewidth]{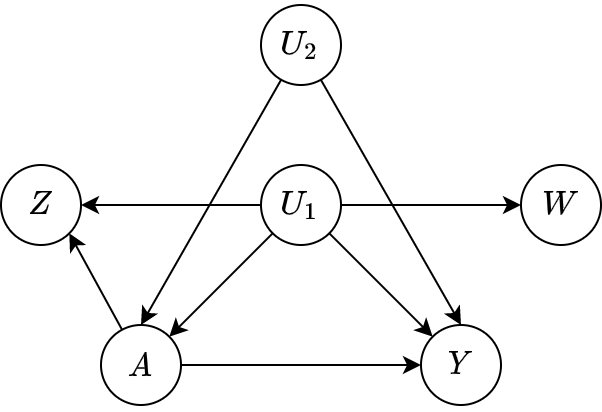}
    \caption{DAG encoding causal relationships among variables in \eqref{eq: U2-AY} in which U-relevance is violated}
    \label{fig: u-relevance}
\end{figure}

In Section 4, we characterize the proximal inference estimator bias in a LSEM under scenarios in which each of $Z$ and $W$ are one-dimensional but $U$ (comprising common causes of any of $A$, $Y$, $Z$, and $W$) has two independent components. We first consider the case where one component of $U$ is an `extra' common cause of $Z$ and $W$ not associated with $A$ or $Y$ (which violates completeness (Assumption 7) and is illustrated in Figure \ref{fig: completeness}), then we consider the case where one component of $U$ is a common cause of $A$ and $Y$ but is not associated with either $Z$ or $W$ (which violates $U$-relevance (Assumption 8) and is illustrated in Figure \ref{fig: u-relevance}). We would argue that it is difficult to guard against violations of Assumptions 7 and 8 arising in this way using subject matter knowledge, making sensitivity analysis for violations of these types particularly valuable.

Additionally, for the settings of Figures \ref{fig: completeness} and \ref{fig: u-relevance}, we compare the bias of the proximal estimator due to violations of Assumptions 7 and 8 to the bias of alternative estimators of the ACE which the analyst might implement under an incorrect unconfoundedness assumption. We consider
\begin{itemize}
    \item [(1)] an outcome regression estimator (referred to as ``OR'') which adjusts for $(Z, W)$ via the g-formula \eqref{eq: g-formula} taking $L = \{Z,W\}$ and specifying outcome regression model $E[Y|A,L]=\beta^T(A,Z,W,AZ,AW)$ 
    \item [(2)] an unadjusted estimator (referred to as ``unadj'') which assumes no unobserved confounding and estimates $E[Y(a)]$ as $\hat{E}[Y|A=a]$ via sample means.
\end{itemize}
\color{black}

In Section 5, we characterize the bias of the proximal estimator in a LSEM where $Z$ and $W$ have the same (arbitrary) number of dimensions and $U$ has at least as many components as either $Z$ or $W$, under the simplifying assumption that the effect of $A$ on $Y$ is not modified by $U$ on the additive scale. This simplifying assumption makes tractable calculations that allow us to develop more general bias formulas for scenarios in which each component of $U$ might have missing arrows into any of $A, Y, Z$, or $W$ in the causal DAG.

\section{Bias Formulas For Two Dimensional $U$}\label{sec: bias-fmlas}
In this section, we will derive formulas for the bias resulting from applying proximal inference under scenarios depicted in Figures 2 and 3 under a LSEM data generating process. 
We base our LSEMs on the data generating process in \cite{miao_identifying_2018}. 


\medskip

Let us consider i.i.d. data generated according to
\begin{equation} \label{eq: U2-AYZW0}
    \begin{split}
        \begin{pmatrix} U \\ X \end{pmatrix}  \sim \mathcal{N} &\left(\begin{pmatrix} 0 \\ 0 \\ 0 \end{pmatrix}, \begin{pmatrix} 1 & 0 & \rho_1 \\ 0 & 1 & \rho_2 \\ \rho_1 & \rho_2 & 1 \end{pmatrix} \right), \ \rho_1, \rho_2 \in \left(-1, 1\right) \\
    \textnormal{logit}(\PP(A = 1 | X,U)) &= \alpha_0 +  \alpha_x X + \boldsymbol{\alpha}_{u}^T U \\
    Z &= \theta_0 + \theta_a A + \theta_x X + \boldsymbol{\theta}_u^T U  + \epsilon_1 \\
    W &= \mu_0 + \mu_x X + \boldsymbol{\mu}_u^T U + \epsilon_2 \\ 
    Y(a) &= \gamma_0 + \gamma_a a + \gamma_x X + \boldsymbol{\gamma}_{u}^T U + \gamma_{au1} aU_1 + \epsilon_3 \\
    \epsilon_1, \epsilon_2, \epsilon_3 &\sim \mathcal{N}(0, 1)
    \end{split}
\end{equation}

\begin{figure}[!h]
    \centering
    \includegraphics[width=0.57\linewidth]{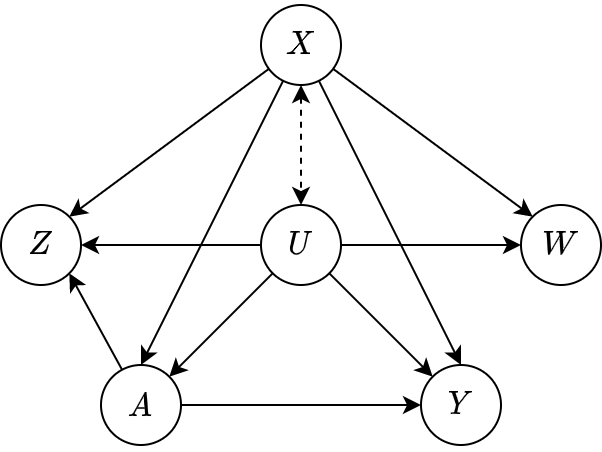}
    \caption{DAG encoding the causal relationships among variables in \eqref{eq: U2-AYZW0}}
    \label{fig: general-figure}
\end{figure}

The causal DAG corresponding to this data generating process can be seen in Figure \ref{fig: general-figure}. Parameter $\boldsymbol{\alpha}_u = \begin{pmatrix} \alpha_{u1} & \alpha_{u2} \end{pmatrix}^T$ encodes the magnitude of confounding, while $\boldsymbol{\theta}_u = \begin{pmatrix} \theta_{u1} & \theta_{u2} \end{pmatrix}^T$ and $\boldsymbol{\mu}_u = \begin{pmatrix} \mu_{u1} & \mu_{u2} \end{pmatrix}^T$ encode the association between confounder $U$ and the NCE/NCO, respectively. We will explore the sensitivity of the proximal inference bias to particular values of $\left(\boldsymbol{\alpha}_u, \boldsymbol{\theta}_u, \boldsymbol{\mu}_u \right)$.

The NCE $Z$ is a post-treatment variable in this DGP. We note that DAGs other than Figure \ref{fig: general-figure} might also be compatible with proximal inference assumptions (e.g., $Z \to A$ or no arrow between $A$ and $Z$, in the absence of other changes). More examples of DAGs compatible with proximal inference assumptions can be found in \cite{shi_selective_2020}. 

\subsection{Base Case: No Violated Assumptions.}

As a sample application of the proximal identification method, we identify the confounding outcome bridge function $h(W, A, X)$ corresponding to the baseline case of one-dimensional $U$ (that is, $\alpha_{u2} = \theta_{u2} = \mu_{u2} = \gamma_{u2} = 0$). For simplicity, we drop the index denoting the first component of $U$. We have the following DGP:
\begin{equation}
    \begin{split}
        \begin{pmatrix} U \\ X \end{pmatrix} &\sim \mathcal{N}\left(\begin{pmatrix} 0 \\ 0 \end{pmatrix}, \begin{pmatrix} 1 & \rho \\ \rho & 1 \end{pmatrix} \right), \ \rho \in \left(-1, 1 \right)\\
    \textnormal{logit}(\PP(A = 1 | X,U)) &= \alpha_0 +  \alpha_x X + \alpha_u U \\
    Z &= \theta_0 + \theta_a A + \theta_x X + \theta_u U  + \epsilon_1 \\
    W &= \mu_0 + \mu_x X + \mu_u U + \epsilon_2 \\ 
    Y(a) &= \gamma_0 + \gamma_a a + \gamma_x X + \gamma_u U + \gamma_{au} aU + \epsilon_3 \\
    \epsilon_1, \epsilon_2, \epsilon_3 &\sim \mathcal{N}(0, 1)
    \end{split}
\end{equation}

The confounding bridge functions that solve both the outcome and treatment bridge function equations take the form
\begin{align}
    h(W, A, X; b) &= b_0 + b_a A + b_w W + b_x X + b_{ax} AX + b_{aw} AW \label{eq: u-outcome-soln}\\
    q(Z, A, X; t) &= 1 + \exp\left\{(-1)^{1-A} (t_0 + t_z Z + t_a A + t_x X ) \right\} \label{eq: u-treatment-soln}
\end{align}
with fitted parameters
\small
\begin{gather}
    (b_0, b_a, b_x, b_w, b_{ax}, b_{aw}) = \left(\gamma_0 - \frac{\mu_0 \gamma_u}{\mu_u}, \gamma_a - \frac{\mu_0 \gamma_{au}}{\mu_u}, \gamma_x - \frac{\mu_x \gamma_u}{\mu_u}, \frac{\gamma_u}{\mu_u}, -\frac{\mu_x \gamma_{au}}{\mu_u}, \frac{\gamma_{au}}{\mu_u} \right) \label{eq: u-outcome-params}\\
    (t_0, t_a, t_x, t_z) = \left(-\alpha_0 + \frac{\theta_0}{\theta_u} \alpha_u + \frac{0.5}{\theta_u^2} \alpha_u^2, -\frac{1}{\theta_u^2} \alpha_u^2 + \frac{\theta_a}{\theta_u} \alpha_u, \frac{\theta_x}{\theta_u} \alpha_u - \alpha_x, -\frac{\alpha_u}{\theta_u} \right) \label{eq: u-treatment-params}
\end{gather}
\normalsize

The proof for the correctness and uniqueness of the above bridge functions is provided in Appendix \hyperref[unobs-u-bridge]{\ref{unobs-u-bridge}}. The proximal g-formula using either bridge function yields an unbiased estimate of the ACE.

\subsection{Violations of Proximal Inference Assumptions}

We examine setup \eqref{eq: U2-AYZW0} for two-dimensional $U$, which implies that at least one of vectors $(\alpha_{u2},\gamma_{u2})$ and $(\theta_{u2}, \mu_{u2})$ has all nonzero entries.

In the case where $\boldsymbol{\theta}_u$ is nonzero (that is, there is a nonzero association between the NCE $Z$ and at least one component of $U$), the following theorem holds:
\begin{theorem}\label{theorem: Completeness violation}
If $\boldsymbol{\theta}_u$ is nonzero (i.e., $Z$ is $U$-relevant), then the LSEM \eqref{eq: U2-AYZW0} with Gaussian $(X,U)$ violates completeness assumption \ref{assumption: Completeness}(a).
\end{theorem}

A proof of Theorem \ref{theorem: Completeness violation} which constructs a counterexample function $g(U)$ for assumption \ref{assumption: Completeness}(a) is provided in Appendix \ref{sec: U2-ZW Completeness pf}. By Theorem \ref{theorem: outcome-bridge}, we know that violating assumption \ref{assumption: Completeness}(a) leads to a potentially biased ACE estimate as the outcome confounding bridge function $\hat h(W, A, X)$ resulting from the GMM procedure no longer satisfies \eqref{eq: outcome-bridge-U}. In the upcoming sections, we will derive formulas for the resulting bias in the above LSEM when $Z$ is $U$-relevant, for the particular cases
\begin{itemize}
    \item $\alpha_{u2} = \gamma_{u2} = 0$ and $\theta_{u2}, \mu_{u2} \ne 0$  (section \ref{section: association btw negative controls}) 
    \item $\theta_{u2} = \mu_{u2} = 0$ and $\alpha_{u2}, \gamma_{u2} \ne 0$ (section \ref{section: partial U-relevance})
\end{itemize}
The two cases were treated separately for simplicity, but they may be combined into a general sensitivity analysis in the context where either vector $(\alpha_{u2}, \gamma_{u2})$ or $(\theta_{u2}, \mu_{u2})$ has all nonzero entries (and the two cases are not mutually exclusive). 

\subsubsection{Completeness Violation: Association between Negative Controls through $U = (U_1, U_2)$ (as in Figure \ref{fig: completeness})}\label{section: association btw negative controls}

For simplicity, we exclude $X$ from these computations. Let us consider i.i.d. data generated according to
\begin{equation} \label{eq: U2-ZW}
    \begin{split}
        U &\sim \mathcal{N}\left(\begin{pmatrix} 0 \\ 0 \end{pmatrix}, \begin{pmatrix} 1 & 0 \\ 0 & 1 \end{pmatrix} \right) \\
    \textnormal{logit}(\PP(A = 1 | U)) &= \alpha_0 + \alpha_{u1} U_1 \\
    Z &= \theta_0 + \theta_a A + \boldsymbol{\theta}_{u}^T U  + \epsilon_1 \\
    W &= \mu_0 + \boldsymbol{\mu}_{u}^T U + \epsilon_2 \\ 
    Y(a) &= \gamma_0 + \gamma_a a + \gamma_{u1} U_1 + \gamma_{au1} aU_1 + 2\epsilon_3 \\
    \epsilon_1, \epsilon_2, \epsilon_3 &\sim \mathcal{N}(0, 1)
    \end{split}
\end{equation}
where $\boldsymbol{\theta}_u, \boldsymbol{\mu}_u$ have all non-zero entries. 

\color{black}

From Theorem \ref{theorem: Completeness violation}, we know that the above setup satisfies all assumptions except \ref{assumption: Completeness}(a). Thus, proceeding to solve for the parameters $b$ of a linear outcome bridge function (which is the functional form an investigator who was unaware of $U_2$ would select) will lead to a biased estimate of the average treatment effect, even if the linear bridge function is correctly specified. The following theorem (see appendix \ref{sec: U2-ZW Completeness pf} for a proof) provides a formula for this bias under a linear outcome bridge function:
\begin{theorem}
If $(A,Y) \indep U_2 \ | \ U_1$, then fitting a linear outcome bridge function $h(W,A,X) = b_0 + b_a A + b_w W + b_x X + b_{ax} AX + b_{aw} AW$ under LSEM \eqref{eq: U2-AYZW0} yields a proximal outcome estimator bias equal to 
\begin{equation}
\begin{split}
    \delta_{POR} =& \frac{\E[AU_1]}{\E[A](1-\E[A])} \frac{\theta_{u2}}{\theta_{u1}} \mu_{u2} \cdot \\
    &\cdot \left[\frac{(1-\E[A]) S_2}{\mu_{u1} + S_2 \cdot \frac{\theta_{u2}}{\theta_{u1}}\mu_{u2}} \gamma_{au1} + \left(\frac{\E[A] S_1}{\mu_{u1} + S_1 \cdot \frac{\theta_{u2}}{\theta_{u1}}\mu_{u2}} + \frac{(1-\E[A]) S_2}{\mu_{u1} + S_2 \cdot \frac{\theta_{u2}}{\theta_{u1}}\mu_{u2}} \right) \gamma_{u1} \right]
    \end{split}
\end{equation}
\normalsize
where
\begin{align*}
    S_1 &= \frac{(1-\E[A])^2}{(1-\E[A])(1-\E[AU_1^2])-\E[AU_1]^2}\\
    S_2 &= \frac{\E[A]^2}{\E[A]\E[AU_1^2]-\E[AU_1]^2}
\end{align*}
In particular, for $\gamma_{au1} = 0$, the bias can be written as 
\begin{equation}
    \delta_{POR} = \frac{\E[AU_1]}{\E[A](1-\E[A])} \frac{\theta_{u2}}{\theta_{u1}} \mu_{u2} \cdot \left(\frac{\E[A] S_1}{\mu_{u1} + S_1 \frac{\theta_{u2}}{\theta_{u1}}\mu_{u2}} + \frac{(1-\E[A]) S_2}{\mu_{u1} + S_2 \cdot \frac{\theta_{u2}}{\theta_{u1}}\mu_{u2}} \right) \gamma_{u1} 
\end{equation}
\end{theorem}

In the remaining theoretical analysis of this case, we make the additional simplifying assumption $\gamma_{au1} = 0$. The general case $\gamma_{au1} \ne 0$ will be considered in the numerical experiments of Section \ref{sec: numerical-experiments}, but we restrict ourselves here for clarity.

\paragraph{Under setup $\gamma_{au1} = 0$: }

By comparison, the bias resulting from the non-proximal g-computation estimator regressing $Y$ onto $(1,Z,W,A,AZ,AW)$ is
\begin{equation}
    \begin{split}
        \delta_{OR} = \frac{\left(\frac{\E[AU_1]}{\E[A]}-\frac{1-\E[A] }{S_2} \theta_a\theta_{u1} \right)\left(1+\mu_{u2}^2 \right) + \left(\frac{\E[AU_1]}{\E[A]}-\frac{1-\E[A] }{S_2} \theta_a\frac{\mu_{u1} \mu_{u2}}{\theta_{u2}} \right) \theta_{u2}^2}{\left(1+\frac{1}{S_2} \theta_{u1}^2 \right)\left(1+\mu_{u2}^2 \right) + \left(1+\frac{1}{S_2}\mu_{u1}^2 \right)\left(1+\theta_{u2}^2 \right) - \left(1 + \frac{2}{S_2} \theta_{u1} \mu_{u1} \theta_{u2} \mu_{u2}\right)} \gamma_{u1}
    \end{split}
\end{equation}
Additionally, the bias resulting from regressing $Y$ onto $A$ is
\begin{equation}
    \delta_{unadj} = \frac{\E[AU_1]}{\E[A] (1 - \E[A])} \gamma_{u1}
\end{equation}
The proofs for the non-proximal g-computation biases can be found in Appendix \ref{sec: g-bias-1}.

It turns out that, under certain setups $\left(\theta_u, \mu_u \right)$ denoting the strengths of association between $(Z,W)$ and $U$, we are guaranteed to obtain less bias from the proximal g-computation estimator than from the unadjusted regression estimator. We formalize these setups in the following theorem:
\begin{theorem}\label{theorem: bias-ZW-comparison}
Under setup $\gamma_{au1} = 0$, the proximal g-computation bias $\delta_{POR}$ and the unadjusted esimator bias $ \delta_{unadj}$ can be compared as follows:
\begin{itemize}
    \item [(i)] If $\theta_{u1} \mu_{u1}$ and $\theta_{u2} \mu_{u2}$ have the same sign (both positive or both negative), then $\left|\delta_{POR} \right| < \left|\delta_{unadj} \right|$.
    \item [(ii)] If $\theta_{u1} \mu_{u1}$ and $\theta_{u2} \mu_{u2}$ have different signs, then \\ $\begin{cases} \left| \delta_{POR} \right| > \left|\delta_{unadj} \right| &\textnormal{if $\frac{\theta_{u1}\mu_{u1}}{\theta_{u2}\mu_{u2}} > -S_1 (1-\E[A]) - S_2 \E[A]$} \\
    \left|\delta_{POR} \right| < \left|\hat \delta_{unadj} \right| &\textnormal{if $\frac{\theta_{u1}\mu_{u1}}{\theta_{u2}\mu_{u2}} < -S_1 (1-\E[A]) - S_2 \E[A]$} \end{cases}$
\end{itemize}
\end{theorem}
The proof of Theorem \ref{theorem: bias-ZW-comparison} can be found in Appendix \ref{sec: bias-comparisons-ZW}.

\subsubsection{Partial $U$-relevance for Two-Dimensional Unobserved Confounder $U$ (as in Figure 3)}\label{section: partial U-relevance}

For simplicity, we exclude $X$ from subsequent computations. Let us consider i.i.d. data generated according to
\begin{equation} \label{eq: U2-AY}
    \begin{split}
        U &\sim \mathcal{N}\left(\begin{pmatrix} 0 \\ 0 \end{pmatrix}, \begin{pmatrix} 1 & 0 \\ 0 & 1 \end{pmatrix} \right) \\  
    \textnormal{logit}(\PP(A = 1 | U)) &= \alpha_0 +  \boldsymbol{\alpha}_{u}^T U \\
    Z &= \theta_0 + \theta_a A + \theta_{u1} U_1  + \epsilon_1 \\
    W &= \mu_0 + \mu_{u1} U_1 + \epsilon_2 \\ 
    Y(a) &= \gamma_0 + \gamma_a a + \boldsymbol{\gamma}_{u}^T U + \gamma_{au1} aU_1 + \epsilon_3 \\
    \epsilon_1, \epsilon_2, \epsilon_3 &\sim \mathcal{N}(0, 1)
    \end{split}
\end{equation}
where $\boldsymbol{\alpha}_u, \boldsymbol{\gamma}_u$ have all non-zero entries. 

From Theorem \ref{theorem: Completeness violation}, we know that the above setup violates assumption \ref{assumption: Completeness}(a). In addition, we do not have a proof identifying the true outcome confounding bridge function, so fitting a linear model might also be misspecified. The following theorem (see Appendix \ref{sec: bias UW-AY pf} for a proof) provides a formula for this bias under a linear bridge function, which can be used in sensitivity analysis. 

\begin{theorem} \label{proof: bias-V-AY}
If $(Z, W) \indep U_2 \ | \ (A, U_1)$, then fitting a linear outcome bridge function $h(W,A,X) = b_0 + b_a A + b_w W + b_x X + b_{ax} AX + b_{aw} AW$ under LSEM \eqref{eq: U2-AYZW0} yields a proximal outcome estimator bias equal to
\begin{equation}
    \delta_{POR} = - \frac{\left(\E[AU_1^2]\left(1-\E[AU_1^2]\right) - \E[AU_1]^2\right) \E[AU_2]}{\left(\E[A]\E[AU_1^2]-\E[AU_1]^2\right)\left(\left(1-\E[A]\right)\left(1-\E[AU_1^2]\right)-\E[AU_1]^2\right)} \gamma_{u2}
\end{equation}
\end{theorem}

By comparison, the bias resulting from the non-proximal g-computation estimator regressing $Y$ onto $(1,Z,W,A,AZ,AW)$ can be obtained as in Section \ref{section: association btw negative controls}, but we omit the formula here due to space constraints and only include this estimate in numerical experiments. 

If we are not considering the proximal estimator, then we may not consider adjusting for the post-exposure variables $Z$ and $W$. In this case, the bias resulting from regressing $Y$ onto $A$ is
\begin{equation}
    \delta_{unadj} = \frac{\E[AU_1]}{\E[A]\left(1-\E[A] \right)} \gamma_{u1} + \frac{\E[AU_1]}{1-\E[A]} \gamma_{au1} + \frac{\E[AU_2]}{\E[A]\left(1-\E[A] \right)} \gamma_{u2}
\end{equation}
\begin{remark}
We note that, while the proximal estimator bias formula depends only on $\gamma_{u2}$, $\alpha_{u2}$ (through $\E[AU_2]$), and $\alpha_{u1}$ (through $\E[AU_1]$), the unadjusted estimator bias depends additionally on parameters $\gamma_{u1}$ and $\gamma_{au1}$ governing the strength of confounding introduced by $U_1$ in the non-proximal case.  
\end{remark}
\color{black}

\section{Bias Formulas in Arbitrary Dimension with No Confounder-Treatment Interaction}

We additionally look into a simplified case where $\gamma_{au} = 0$ -- that is, the confounder is not an effect modifier. Moreover, we assume that the analyst is aware of the lack of interaction between $A$ and $U$ in the true outcome model, so we consider a simplified bridge function model $h(W,A,X) = b_0 + b_a A + b_w^T W + b_x^T X$ as input. This assumption allows us to more easily obtain bias formulas in the general case of multi-dimensional $Z,W,U,X$ with $\left(dim(Z), dim(W), dim(U), dim(X) \right) = (m, n, p, q)$, for certain relationships between $m, n, p, q$. For simplicity, we assume that the unobserved and observed confounders $(U, X)$ jointly follow a multivariate normal distribution with mean $\boldsymbol{0}_{p+q}$, $Var(U_i) = Var(X_j) = 1$ for all $i = 1,\ldots,p$, $j = 1,\ldots,q$ (under a potential transformation), and some appropriate PSD covariance matrix such that $Cov(U, X) = \rho \in (-1, 1)^{p \times q}$.

Let us consider i.i.d. data generated according to 
\begin{equation} \label{eq: U2-AYZW-complete}
    \begin{split}
        \begin{pmatrix} U \\ X \end{pmatrix} &\sim \mathcal{N}\left(\begin{pmatrix} \boldsymbol{0}_p \\ \boldsymbol{0}_q \end{pmatrix}, \begin{pmatrix} \boldsymbol{I}_p & \rho \\ \rho^T & \boldsymbol{\Sigma}_x \end{pmatrix} \right), \ \rho \in \left(-1, 1 \right)^{p \times q}\\
    \textnormal{logit}(\PP(A = 1 | U, X)) &= \alpha_0 +  \alpha_{u}^T U + \alpha_{x}^T X \\
    Z &= \theta_0 + \theta_a A + \theta_u^T U  + \theta_x^T X + \epsilon_1 \\
    W &= \mu_0 + \mu_u^T U + \mu_x^T X + \epsilon_2 \\ 
    Y(a) &= \gamma_0 + \gamma_a a +  \gamma_{u}^T U + \gamma_x X + \epsilon_3 \\
    \epsilon_1, \epsilon_2, \epsilon_3 &\sim \mathcal{N}(0, 1)
    \end{split}
\end{equation}

The following theorem provides a formula for the proximal identification bias under a linear bridge function.
\begin{theorem} \label{theorem: general-bias}
Let $\E[AU] = \left(\E[AU_1],\ldots,\E[AU_p] \right)$, $\E[AX] = \left(\E[AX_1],\ldots,\E[AX_p] \right)$, and 
\begin{equation*}
    B = \left(I_p - \rho \Sigma_x^{-1} \rho^T - \frac{\left(\E[AU] - \rho \Sigma_x^{-1} \E[AX] \right) \left(\E[AU]^T - \E[AX]^T \Sigma_x^{-1} \rho^T \right)}{\E[A]\left(1-\E[A]\right) - \E[AX]^T \Sigma_x^{-1} \E[AX]} \right) \theta_u.
\end{equation*} 
If $m = n < p$ and matrix $B^T \mu_u \in \R^{m \times m}$ has full rank, then fitting a linear outcome bridge function $h(W,A,X) = b_0 + b_a A + b_w^T W + b_x^T X$ under LSEM \eqref{eq: U2-AYZW-complete} yields a proximal outcome estimator bias equal to
\begin{equation}\label{sensitivity_bias}
    \delta = \frac{\E[AU]^T - \E[AX]^T \Sigma_x^{-1} \rho^T}{\E[A] \left(1 - \E[A] \right) - \E[AX]^T \Sigma_x^{-1} \E[AX]} \left(I_p - \mu_u \left(B^T \mu_u \right)^{-1} B^T \right) \gamma_u 
\end{equation}
\end{theorem}

A proof of Theorem \ref{theorem: general-bias} can be found in Appendix \ref{sec: general-pi-bias}. 

\begin{remark}\label{remark2} If $m = n = p$ and $B^T \mu_u$ has full rank, then $\delta = 0$. If $p < m$ or $p < n$, then we have a similar discussion as in \cite{shi_multiply_2019-1} where we can either consider the Moore-Penrose inverse of $B^T \mu_u$, or reduce the dimensions of $Z$ and $W$ until they match the dimension of $U$.
\end{remark}

\begin{remark}\label{remark3} Theorem \ref{theorem: general-bias} enables sensitivity analysis. Note that the terms $\E[A]$ and $\E[AX]$ in (\ref{sensitivity_bias}) are straightforwardly estimated from data. Thus, to perform a sensitivity analysis using the bias formula (\ref{sensitivity_bias}), it remains for the analyst to specify the parameters $\E[AU]$ (which is determined by $\alpha_u$), $\mu_u$, $\gamma_u$, and $\rho$. An analyst could specify a distribution over these parameters, which, via (\ref{sensitivity_bias}), would imply a distribution over $\delta$ as each realization of the parameters drawn from the distribution would correspond to a different bias $\delta$. The range of magnitudes of the parameters governing the strength of association between $U$ and the other variables might be chosen based on the range of magnitudes of associations between the observed variables. And probabilities of zero components in $\theta_u$ and $\mu_u$ could determine the proportion of components of $U$ that contribute to bias from $U-$irrelevance. 
\end{remark}

\section{Numerical Experiments}\label{sec: numerical-experiments}

In this section, we provide numerical examples to illustrate how the bias of different estimators (proximal and non-proximal) varies with the strength and direction of associations between unobserved $U$ and $(A,Y,Z,W)$. Due to the relatively large number of parameters involved in the bias formulas, we fix values $\alpha_0 = \gamma_0 = \theta_0 = \mu_0 = 0$ , $\theta_a = \theta_{u1} = 1$, $\mu_{u1} = 0.5$, $\gamma_{u1} = 1$, $\gamma_{au1} = 1.5$ similar to the simulation DGP in \cite{miao_confounding_2020}. We then analyze bias sensitivity to different values of $(\alpha_{u2}, \theta_{u2}, \mu_{u2}, \gamma_{u2})$, which encode how strongly the proximal identification assumptions are violated in the presence of $U_2$. 
The numerical results and plots in this discussion have been outputted in Mathematica.

\subsection{Completeness Violation: Association between Negative Controls through $U = (U_1, U_2)$}

\subsubsection{Comparison of estimators for $\theta_{u1},\theta_{u2},\mu_{u1},\mu_{u2} > 0$}

\begin{figure}[!h]
\centering
\begin{subfigure}{0.49\textwidth}
  \centering
  \includegraphics[width=\linewidth]{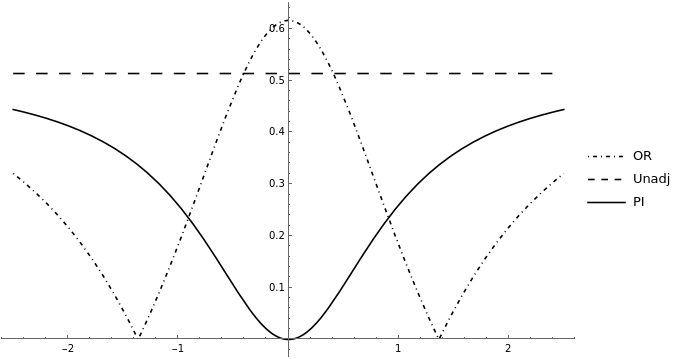}
  \caption{$\mu_{u2} = 0.5 \theta_{u2}$}
  \label{fig:sfig1}
\end{subfigure}
\begin{subfigure}{0.49\textwidth}
  \centering
  \includegraphics[width=\linewidth]{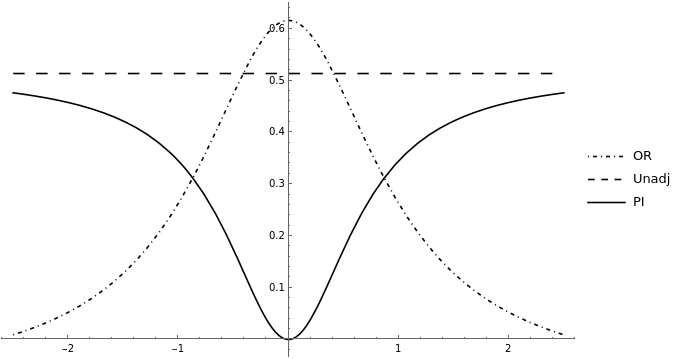}
  \caption{$\mu_{u2} = \theta_{u2}$}
  \label{fig:sfig2}
\end{subfigure}
\begin{subfigure}{0.49\textwidth}
  \centering
  \includegraphics[width=\linewidth]{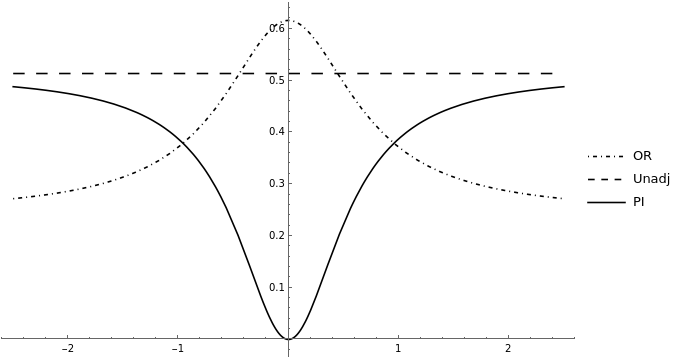}
  \caption{$\mu_{u2} = 1.5 \theta_{u2}$}
  \label{fig:sfig2}
\end{subfigure}
\caption{Plots of the ACE estimate bias for $\alpha_{u1} = 0.3$.}
\label{fig:fig}
\end{figure}

\begin{figure}[!h]
\centering
\begin{subfigure}{0.49\textwidth}
  \centering
  \includegraphics[width=\linewidth]{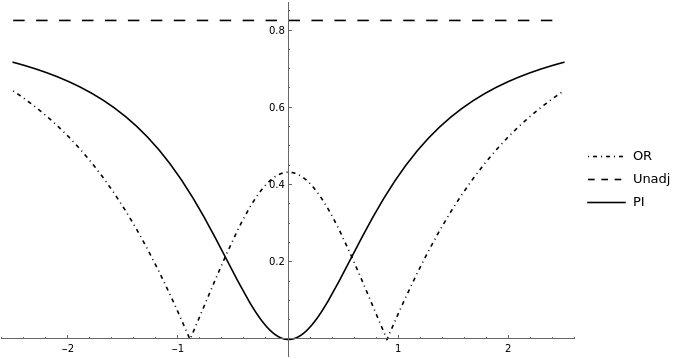}
  \caption{$\mu_{u2} = 0.5 \theta_{u2}$}
  \label{fig:sfig1}
\end{subfigure}
\begin{subfigure}{0.49\textwidth}
  \centering
  \includegraphics[width=\linewidth]{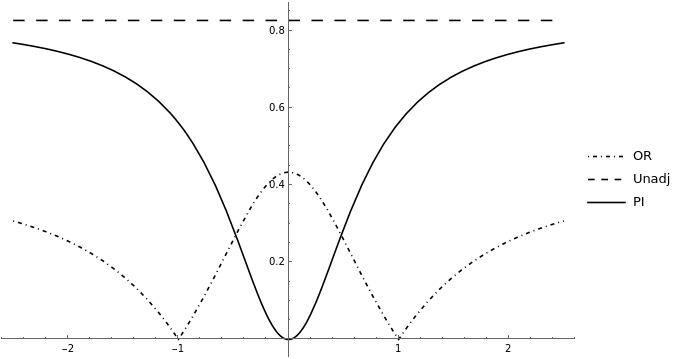}
  \caption{$\mu_{u2} = \theta_{u2}$}
  \label{fig:sfig2}
\end{subfigure}
\begin{subfigure}{0.49\textwidth}
  \centering
  \includegraphics[width=\linewidth]{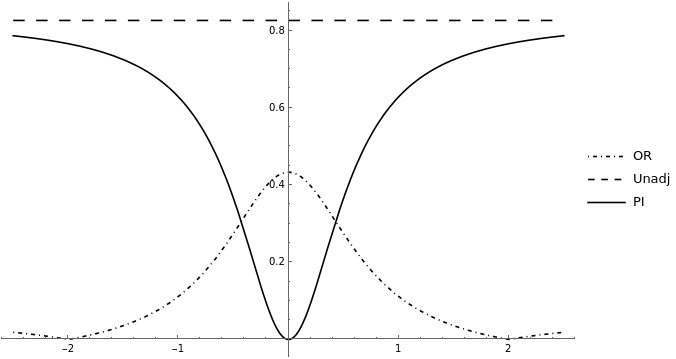}
  \caption{$\mu_{u2} = 1.5 \theta_{u2}$}
  \label{fig:sfig2}
\end{subfigure}
\caption{Plots of the ACE estimate bias for $\alpha_{u1} = 0.5$.}
\label{fig:fig2}
\end{figure}

\begin{figure}[!h]
\centering
\begin{subfigure}{0.49\textwidth}
  \centering
  \includegraphics[width=\linewidth]{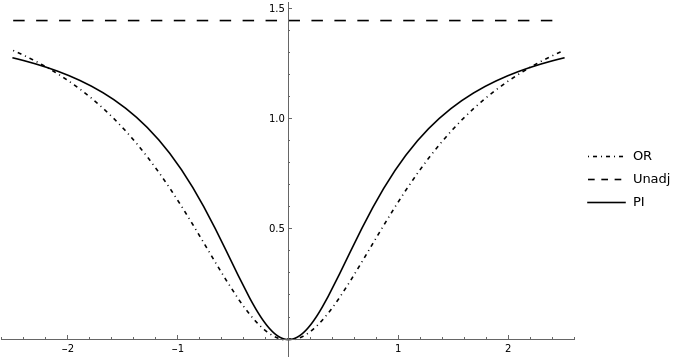}
  \caption{$\mu_{u2} = 0.5 \theta_{u2}$}
  \label{fig:sfig1}
\end{subfigure}
\begin{subfigure}{0.49\textwidth}
  \centering
  \includegraphics[width=\linewidth]{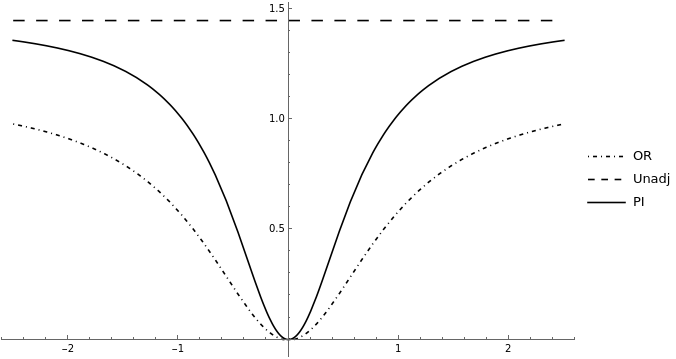}
  \caption{$\mu_{u2} = \theta_{u2}$}
  \label{fig:sfig2}
\end{subfigure}\\
\begin{subfigure}{0.49\textwidth}
  \centering
  \includegraphics[width=\linewidth]{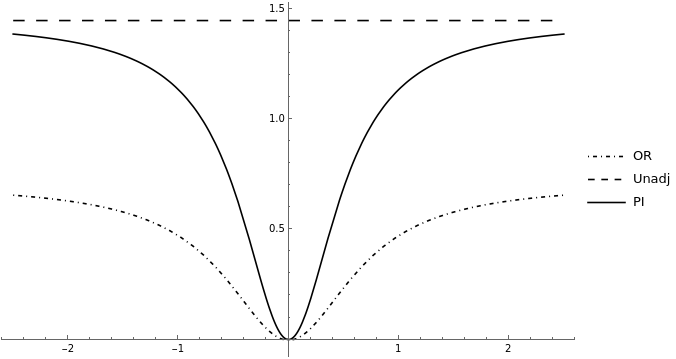}
  \caption{$\mu_{u2} = 1.5 \theta_{u2}$}
  \label{fig:sfig2}
\end{subfigure}
\caption{Plots of the ACE estimate bias for $\alpha_{u1} = 1$.}
\label{fig:fig3}
\end{figure}

\begin{figure}[!h]
\centering
\begin{subfigure}{0.49\textwidth}
  \centering
  \includegraphics[width=\linewidth]{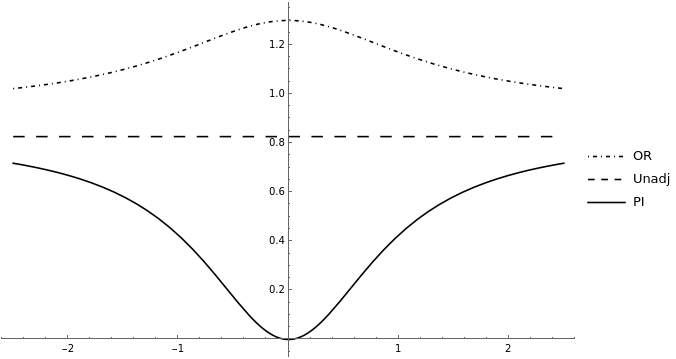}
  \caption{$\mu_{u2} = 0.5 \theta_{u2}$}
  \label{fig:sfig1}
\end{subfigure}
\begin{subfigure}{0.49\textwidth}
  \centering
  \includegraphics[width=\linewidth]{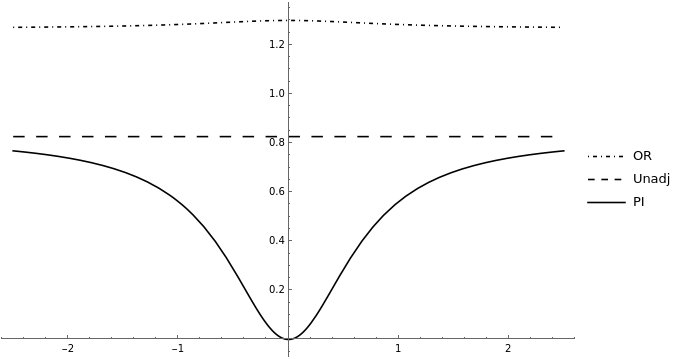}
  \caption{$\mu_{u2} = \theta_{u2}$}
  \label{fig:sfig2}
\end{subfigure}
\begin{subfigure}{0.49\textwidth}
  \centering
  \includegraphics[width=\linewidth]{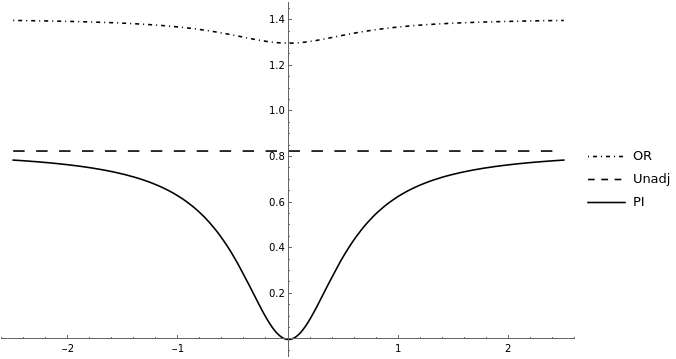}
  \caption{$\mu_{u2} = 1.5 \theta_{u2}$}
  \label{fig:sfig2}
\end{subfigure}
\caption{Plots of the ACE estimate bias for $\alpha_{u1} = -0.3$.}
\label{fig:fig4}
\end{figure}

This is the case considered in Theorem \ref{theorem: bias-ZW-comparison}(i), where it is shown that the bias of the unadjusted estimator always exceeds that of the proximal estimator. Figures \ref{fig:fig}-\ref{fig:fig4} illustrate the change in absolute bias for each of the three estimators. In all figures, we use the following notation:
\begin{itemize}
    \item [(1)] The solid black curve (``PI'') corresponds to the (absolute) proximal estimator bias
    \item [(2)] The dashed curve (``Unadj'') corresponds to the (absolute) unadjusted estimator bias from regressing $Y$ on $A$
    \item [(3)] The dot-dashed curve (``OR'') corresponds to the (absolute) adjusted estimator bias from regressing $Y$ on $(1,Z,W,A,AZ,AW)$
\end{itemize}


Consistent with Theorem 4 (i), the proximal estimator always outperforms the unadjusted estimator. However, we also note that there exist settings $(\alpha_{u1}, \theta_{u2}, \mu_{u2})$ for which the non-proximal adjusted estimate is (significantly) less biased than the proximal one. Although standard criteria for variable adjustment do not usually include post-exposure covariates in the adjustment set \cite{vanderweele_principles_2019}, the strong underlying association between $U$ and the observed proxies $(Z,W)$ might actually help mitigate bias through adjustment.  

\medskip
\subsubsection{Comparison of Estimators for Different Directions of Association Products $\theta_{u1} \mu_{u1}$, $\theta_{u2} \mu_{u2}$}

\paragraph{Same sign of $\theta_{u1} \mu_{u1}$ and $\theta_{u2} \mu_{u2}$:}

\begin{figure}[!h]
\centering
\begin{subfigure}{0.49\textwidth}
  \centering
  \includegraphics[width=\linewidth]{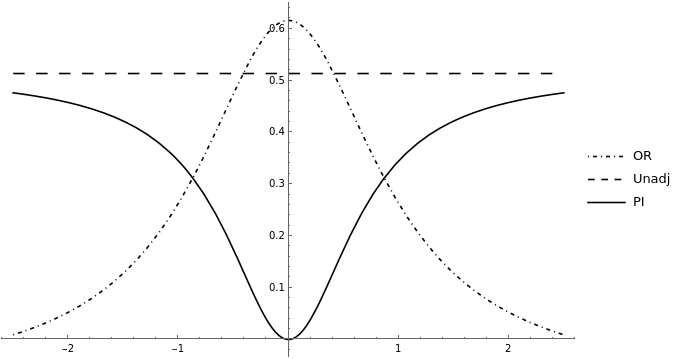}
  \caption{$\alpha_{u1} = 0.3$, $\mu_{u2} = \theta_{u2}$}
  \label{fig:sfig1}
\end{subfigure}
\begin{subfigure}{0.49\textwidth}
  \centering
  \includegraphics[width=\linewidth]{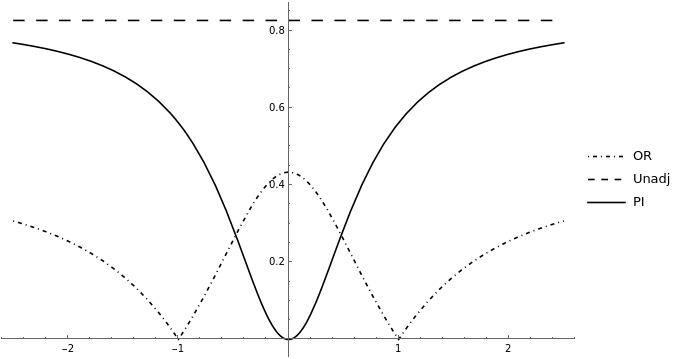}
  \caption{$\alpha_{u1} = 0.5$, $\mu_{u2} = \theta_{u2}$}
  \label{fig:sfig2}
\end{subfigure}
\begin{subfigure}{0.49\textwidth}
  \centering
  \includegraphics[width=\linewidth]{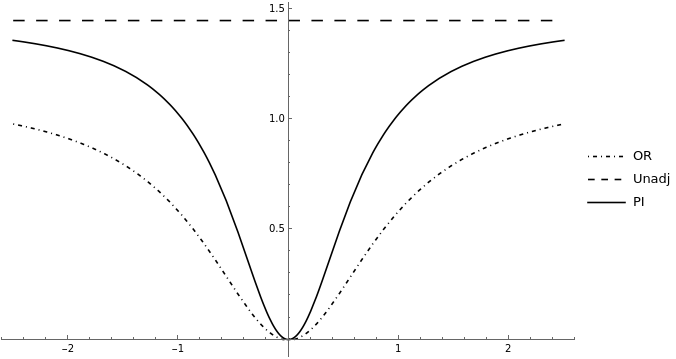}
  \caption{$\alpha_{u1} = 1$, $\mu_{u2} = \theta_{u2}$}
  \label{fig:sfig2}
\end{subfigure}
\begin{subfigure}{0.49\textwidth}
  \centering
  \includegraphics[width=\linewidth]{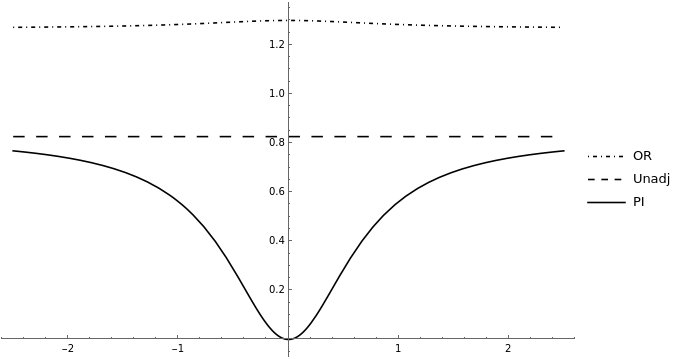}
  \caption{$\alpha_{u1} = -0.5$, $\mu_{u2} = \theta_{u2}$}
  \label{fig:sfig2}
\end{subfigure}
\caption{Plots of the ACE estimate bias for $\theta_{u1} = \mu_{u1} = 0.5$.}
\label{fig:fig5}
\end{figure}

Figure \ref{fig:fig5} illustrates the change in absolute bias for each of the three estimators relative to the value of $\theta_{u1}$, where it is assumed that $\mu_{u2} = \theta_{u2}$ in all cases. We observe that the absolute unadjusted bias is always greater than the proximal estimator bias, which is consistent with the result in Theorem \ref{theorem: bias-ZW-comparison}. Comparisons with the adjusted estimator are not as straightforward, but there seem to exist threshold values of $\alpha_{u1}$ which determine whether the adjusted estimator bias ever exceeds the unadjusted one (such as in Figures \ref{fig:fig5}(a) and (d)). Moreover, when $\alpha_{u1} > 0$, there exists a threshold value of $|\theta_{u2}|$ which determines whether the proximal estimator bias exceeds the adjusted bias.

The parameters used in these plots are: $\alpha_0 = \gamma_0 = \theta_0 = \mu_0 = 0$, $\theta_a = \theta_{u1} = 1$, $\gamma_{u1} = 1$, $\gamma_{ua1} = 1.5$, $\mu_{u1} = 0.5$, $\gamma_a = 0.5$.

\paragraph{Different signs of $\theta_{u1} \mu_{u1}$ and $\theta_{u2} \mu_{u2}$:}

\begin{figure}[!h]
\centering
\begin{subfigure}{0.49\textwidth}
  \centering
  \includegraphics[width=\linewidth]{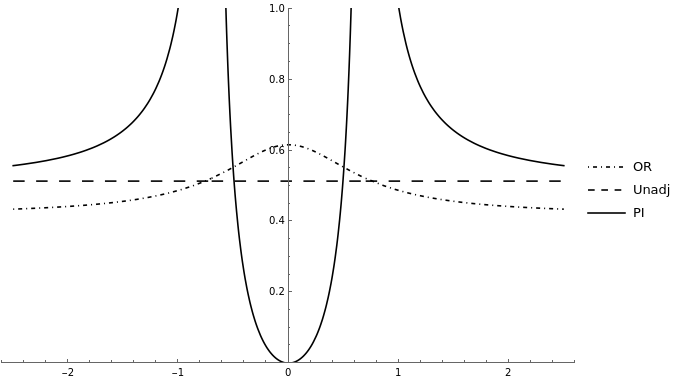}
  \caption{$\alpha_{u1} = 0.3$, $\mu_{u2} = \theta_{u2}$}
  \label{fig:sfig1}
\end{subfigure}
\begin{subfigure}{0.49\textwidth}
  \centering
  \includegraphics[width=\linewidth]{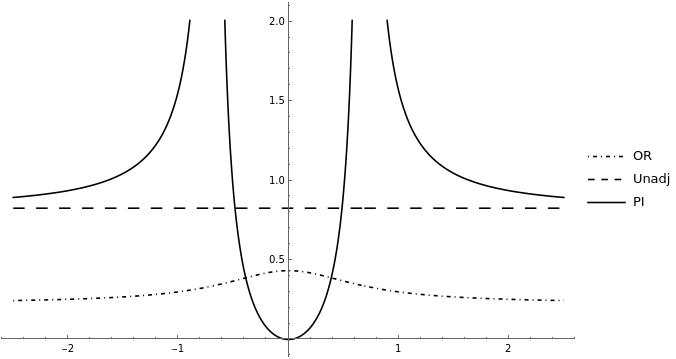}
  \caption{$\alpha_{u1} = 0.5$, $\mu_{u2} = \theta_{u2}$}
  \label{fig:sfig2}
\end{subfigure}
\begin{subfigure}{0.49\textwidth}
  \centering
  \includegraphics[width=\linewidth]{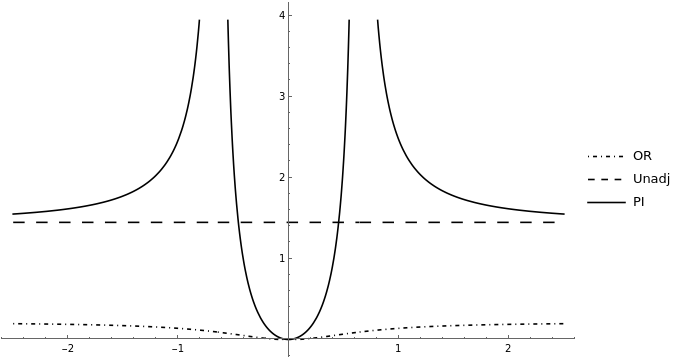}
  \caption{$\alpha_{u1} = 1$, $\mu_{u2} = \theta_{u2}$}
  \label{fig:sfig2}
\end{subfigure}
\begin{subfigure}{0.49\textwidth}
  \centering
  \includegraphics[width=\linewidth]{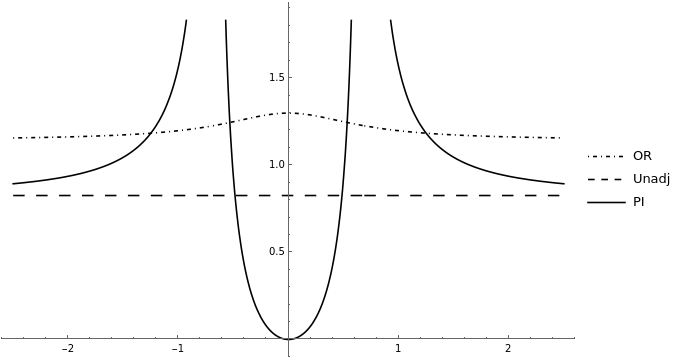}
  \caption{$\alpha_{u1} = -0.5$, $\mu_{u2} = \theta_{u2}$}
  \label{fig:sfig2}
\end{subfigure}
\caption{Plots of the ACE estimate bias for $\theta_{u1} = \mu_{u1} = -0.5$.}
\label{fig:fig6}
\end{figure}

Figure \ref{fig:fig6} illustrates the change in absolute bias for each of the three estimators relative to the value of $\theta_{u1}$, where it is assumed that $\mu_{u2} = \theta_{u2}$ in all cases. We observe that, beyond a certain threshold in the value of $|\theta_{u2}|$, the proximal estimation bias exceeds that of the unadjusted estimator (and even the adjusted estimator bias, for $\alpha_{u1} > 0$), which is consistent with Theorem \ref{theorem: bias-ZW-comparison}.

The parameters used in these plots are: $\alpha_0 = \gamma_0 = \theta_0 = \mu_0 = 0$, $\theta_a = \theta_{u1} = 1$, $\gamma_{u1} = 1$, $\gamma_{ua1} = 1.5$, $\mu_{u1} = -0.5$, $\gamma_a = 0.5$.

\subsection{Partial $U$-relevance for Two-Dimensional Unobserved Confounder $U$}

\paragraph{Same directions of association $\gamma_{u1}, \gamma_{u2}$:}

\begin{figure}[!h]
\centering
\begin{subfigure}{0.49\textwidth}
  \centering
  \includegraphics[width=\linewidth]{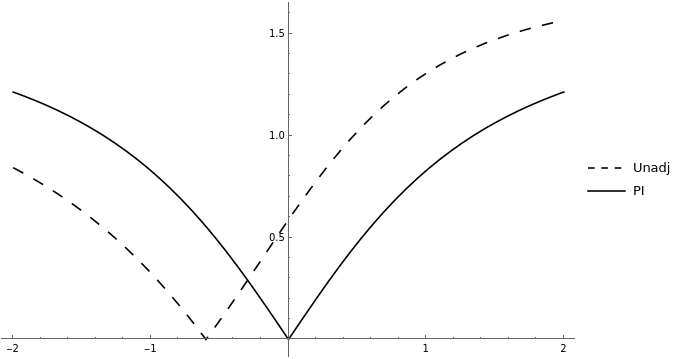}
  \caption{$\alpha_{u1} = 0.3$}
  \label{fig:sfig1}
\end{subfigure}
\begin{subfigure}{0.49\textwidth}
  \centering
  \includegraphics[width=\linewidth]{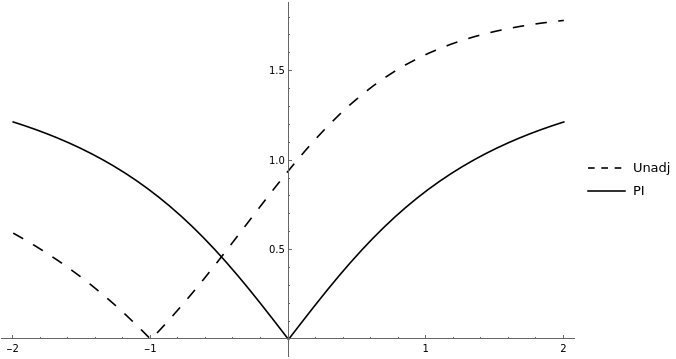}
  \caption{$\alpha_{u1} = 0.5$}
  \label{fig:sfig2}
\end{subfigure}
\begin{subfigure}{0.49\textwidth}
  \centering
  \includegraphics[width=\linewidth]{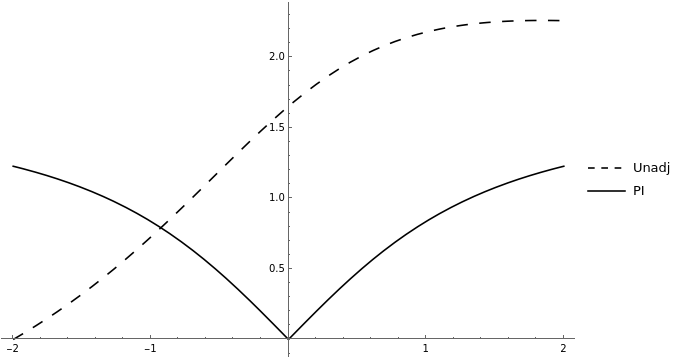}
  \caption{$\alpha_{u1} = 1$}
  \label{fig:sfig3}
\end{subfigure}
\begin{subfigure}{0.49\textwidth}
  \centering
  \includegraphics[width=\linewidth]{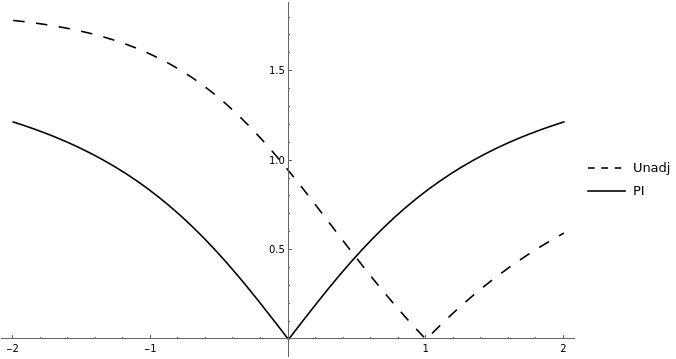}
  \caption{$\alpha_{u1} = -0.5$}
  \label{fig:sfig4}
\end{subfigure}
\caption{Plots of the ACE estimate bias for $\gamma_{u1} = 1.5, \gamma_{u2} > 0$.}
\label{fig:fig7}
\end{figure}

Figure \ref{fig:fig7} illustrates the change in absolute bias for the proximal and unadjusted estimators relative to the value of $\alpha_{u2}$.

The parameters used in these plots are: $\alpha_0 = \gamma_0 = \theta_0 = \mu_0 = 0$, $\theta_a = \theta_{u1} = 1$, $\mu_{u1} = 1$, $\gamma_{au1} = 1$, $\gamma_{u1} = 1.5$, $\gamma_{u2} = 1$, $\gamma_a = 0.5$.

\paragraph{Opposite directions of association $\gamma_{u1}, \gamma_{u2}$:}

\begin{figure}[!h]
\centering
\begin{subfigure}{0.49\textwidth}
  \centering
  \includegraphics[width=\linewidth]{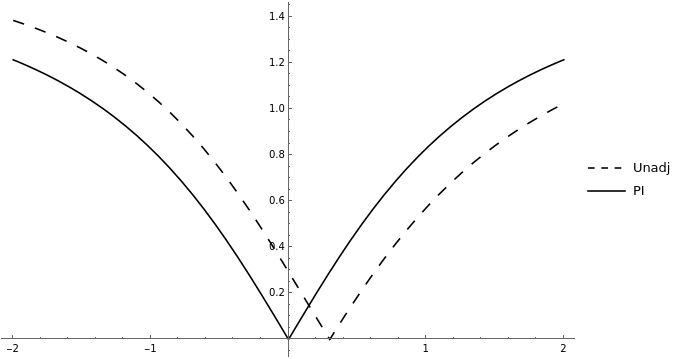}
  \caption{$\alpha_{u1} = 0.3$}
  \label{fig:sfigneg1}
\end{subfigure}
\begin{subfigure}{0.49\textwidth}
  \centering
  \includegraphics[width=\linewidth]{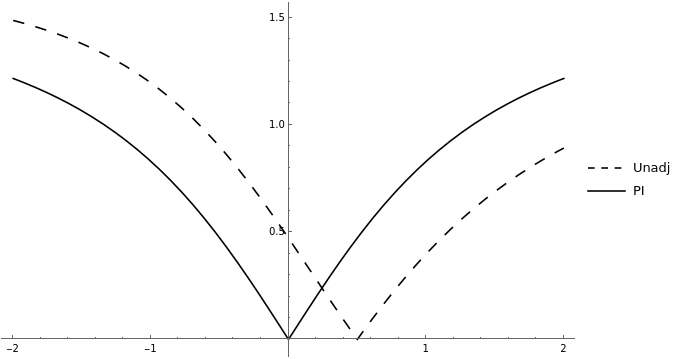}
  \caption{$\alpha_{u1} = 0.5$}
  \label{fig:sfigneg2}
\end{subfigure}
\begin{subfigure}{0.49\textwidth}
  \centering
  \includegraphics[width=\linewidth]{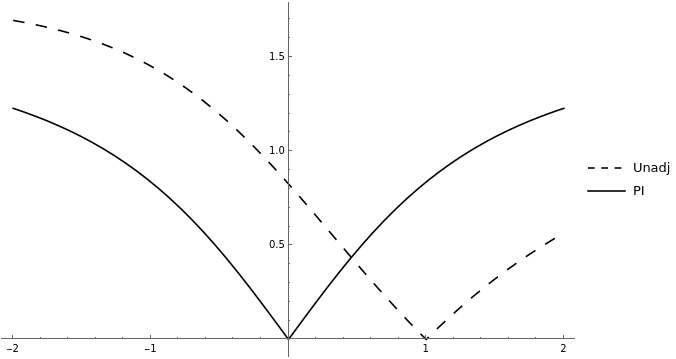}
  \caption{$\alpha_{u1} = 1$}
  \label{fig:sfigneg3}
\end{subfigure}
\begin{subfigure}{0.49\textwidth}
  \centering
  \includegraphics[width=\linewidth]{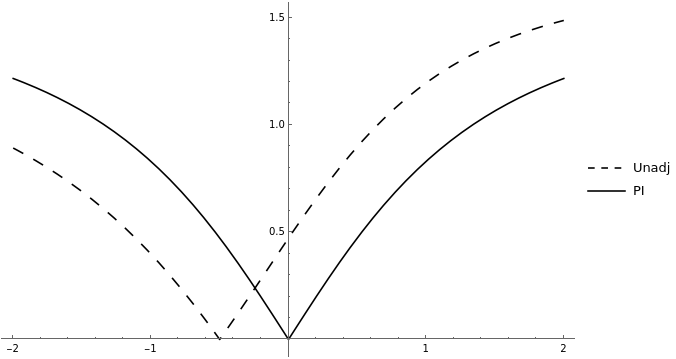}
  \caption{$\alpha_{u1} = -0.5$}
  \label{fig:sfigneg4}
\end{subfigure}
\caption{Plots of the ACE estimate bias for $\gamma_{u1} = -1.5, \gamma_{u2} > 0$.}
\label{fig:fig8}
\end{figure}

Figure \ref{fig:fig8} illustrates the change in absolute bias for the proximal and unadjusted estimators relative to the value of $\alpha_{u2}$.

The parameters used in these plots are: $\alpha_0 = \gamma_0 = \theta_0 = \mu_0 = 0$, $\theta_a = \theta_{u1} = 1$, $\mu_{u1} = 1$, $\gamma_{ua1} = 1$, $\gamma_{u1} = -1.5$, $\gamma_{u2} = 1$, $\gamma_a = 0.5$.

The distributions of bias appear almost shifted by translation. We observe a reversal in which estimator has less bias compared to the case of $\gamma_{u1},\gamma_{u2} > 0$.

\section{Discussion}
By deriving bias formulas for proximal inference estimators under violations of completeness and $U$-relevance under a LSEM, we begin to gain some insight into the sensitivity of the proximal inference estimator to these sources of bias. For example, under some settings, it is possible for completeness violations alone without any failure of $U$-relevance (i.e. too many common causes of the negative control exposure and negative control outcome) to lead the proximal inference estimator to be arbitrarily more biased than an unadjusted estimator completely subject to unobserved confounding (see Figure \ref{fig:fig6}). However, under a LSEM, if the unobserved confounder leads to a positive association between the negative control exposure and negative control outcome, then under a completeness violation with full $U$-relevance the proximal inference estimator is guaranteed to perform better than an unadjusted one, no matter how strong the completeness violation (as shown in Theorem \ref{theorem: bias-ZW-comparison}). A tentative rule of thumb for the design of proximal inference studies, pending additional evidence from other data generating processes, should be to select negative control exposures and outcomes that are positively associated. 

Additionally, as discussed in Remark \ref{remark3}, we can use our bias formula results (in particular, \eqref{sensitivity_bias}) to devise schemes for sensitivity analyses of proximal inference studies. While \eqref{sensitivity_bias} was derived under the strong assumptions that the data generating process is a LSEM and $U$ is not an effect modifier, an analyst might reasonably conduct a sensitivity analysis using (\ref{sensitivity_bias}) as we described even if they did not believe their data were generated by a LSEM (and did not construct their proximal inference estimators under that assumption), and even if they did not believe that $U$ is not an effect modifier. There is a long history of unrealistic simplifying assumptions in sensitivity analysis. For example, \cite{rosenbaum_sensitivity_1987} and \cite{vanderweele_bias_2011} both assume a one-dimensional binary confounder for tractable sensitivity analysis of no unobserved confounding. Later, \cite{ding_sensitivity_2016} developed an approach that made far fewer restrictions. We are in the early stages of proximal inference, so we currently need to settle for sensitivity analyses that make strong simplifying assumptions.

\newpage

\bibliographystyle{siam}
\bibliography{references}

\begin{thebibliography}{10}

\bibitem{brookhart_preference-based_2007}
{\sc M.~A. Brookhart and S.~Schneeweiss}, {\em Preference-based instrumental
  variable methods for the estimation of treatment effects: assessing validity
  and interpreting results}, The International Journal of Biostatistics, 3
  (2007), p.~Article 14.

\bibitem{cornfield_smoking_2009}
{\sc J.~Cornfield, W.~Haenszel, E.~C. Hammond, A.~M. Lilienfeld, M.~B. Shimkin,
  and E.~L. Wynder}, {\em Smoking and lung cancer: recent evidence and a
  discussion of some questions*}, International Journal of Epidemiology, 38
  (2009), pp.~1175--1191.

\bibitem{cui_semiparametric_2020}
{\sc Y.~Cui, H.~Pu, X.~Shi, W.~Miao, and E.~T. Tchetgen}, {\em Semiparametric
  proximal causal inference}, arXiv:2011.08411 [math, stat],  (2020).
\newblock arXiv: 2011.08411.

\bibitem{ding_sensitivity_2016}
{\sc P.~Ding and T.~J. VanderWeele}, {\em Sensitivity {Analysis} {Without}
  {Assumptions}}, Epidemiology (Cambridge, Mass.), 27 (2016), pp.~368--377.

\bibitem{greenland_identifiability_1986}
{\sc S.~Greenland and J.~M. Robins}, {\em Identifiability, exchangeability, and
  epidemiological confounding}, International Journal of Epidemiology, 15
  (1986), pp.~413--419.

\bibitem{lash_good_2014}
{\sc T.~L. Lash, M.~P. Fox, R.~F. MacLehose, G.~Maldonado, L.~C. McCandless,
  and S.~Greenland}, {\em Good practices for quantitative bias analysis},
  International Journal of Epidemiology, 43 (2014), pp.~1969--1985.

\bibitem{miao_identifying_2018}
{\sc W.~Miao, Z.~Geng, and E.~J. Tchetgen~Tchetgen}, {\em Identifying causal
  effects with proxy variables of an unmeasured confounder}, Biometrika, 105
  (2018), pp.~987--993.

\bibitem{miao_confounding_2020}
{\sc W.~Miao, X.~Shi, and E.~T. Tchetgen}, {\em A {Confounding} {Bridge}
  {Approach} for {Double} {Negative} {Control} {Inference} on {Causal}
  {Effects}}, arXiv:1808.04945 [stat],  (2020).
\newblock arXiv: 1808.04945.

\bibitem{rambachan_more_2022}
{\sc A.~Rambachan and J.~Roth}, {\em A {More} {Credible} {Approach} to
  {Parallel} {Trends}},  (2022), p.~79.

\bibitem{robins_sensitivity_2000}
{\sc J.~M. Robins, A.~Rotnitzky, and D.~O. Scharfstein}, {\em Sensitivity
  {Analysis} for {Selection} bias and unmeasured {Confounding} in missing
  {Data} and {Causal} inference models}, in Statistical {Models} in
  {Epidemiology}, the {Environment}, and {Clinical} {Trials}, M.~E. Halloran
  and D.~Berry, eds., The {IMA} {Volumes} in {Mathematics} and its
  {Applications}, New York, NY, 2000, Springer, pp.~1--94.

\bibitem{rosenbaum_sensitivity_1987}
{\sc P.~R. Rosenbaum}, {\em Sensitivity {Analysis} for {Certain} {Permutation}
  {Inferences} in {Matched} {Observational} {Studies}}, Biometrika, 74 (1987),
  pp.~13--26.
\newblock Publisher: [Oxford University Press, Biometrika Trust].

\bibitem{rosenbaum_assessing_1983}
{\sc P.~R. Rosenbaum and D.~B. Rubin}, {\em Assessing {Sensitivity} to an
  {Unobserved} {Binary} {Covariate} in an {Observational} {Study} with {Binary}
  {Outcome}}, Journal of the Royal Statistical Society. Series B
  (Methodological), 45 (1983), pp.~212--218.
\newblock Publisher: [Royal Statistical Society, Wiley].

\bibitem{rubin_estimating_1974}
{\sc D.~B. Rubin}, {\em Estimating causal effects of treatments in randomized
  and nonrandomized studies}, Journal of Educational Psychology, 66 (1974),
  pp.~688--701.
\newblock Place: US Publisher: American Psychological Association.

\bibitem{shi_multiply_2019-1}
{\sc X.~Shi, W.~Miao, J.~C. Nelson, and E.~J.~T. Tchetgen}, {\em Multiply
  {Robust} {Causal} {Inference} with {Double} {Negative} {Control} {Adjustment}
  for {Categorical} {Unmeasured} {Confounding}}, Sept. 2019.
\newblock arXiv:1808.04906 [stat].

\bibitem{shi_selective_2020}
{\sc X.~Shi, W.~Miao, and E.~T. Tchetgen}, {\em A {Selective} {Review} of
  {Negative} {Control} {Methods} in {Epidemiology}}, arXiv:2009.05641 [stat],
  (2020).
\newblock arXiv: 2009.05641.

\bibitem{tchetgen_introduction_2020}
{\sc E.~J.~T. Tchetgen, A.~Ying, Y.~Cui, X.~Shi, and W.~Miao}, {\em An
  {Introduction} to {Proximal} {Causal} {Learning}}, arXiv:2009.10982 [stat],
  (2020).
\newblock arXiv: 2009.10982.

\bibitem{vanderweele_principles_2019}
{\sc T.~J. VanderWeele}, {\em Principles of confounder selection}, European
  Journal of Epidemiology, 34 (2019), pp.~211--219.

\bibitem{vanderweele_bias_2011}
{\sc T.~J. Vanderweele and O.~A. Arah}, {\em Bias formulas for sensitivity
  analysis of unmeasured confounding for general outcomes, treatments, and
  confounders}, Epidemiology (Cambridge, Mass.), 22 (2011), pp.~42--52.

\end{thebibliography}

\newpage

\appendix

\appendixhead

\section{Bridge Function Parameters for Post-Treatment NCE}

\subsection{Bridge functions derivation for one-dimensional unobserved $U$ -- case of no violations} \label{unobs-u-bridge}

We identify coefficients $(b_0, b_a, b_x, b_w, b_{ax}, b_{aw})$ and $(t_0, t_a, t_x, t_z)$ such that
\begin{align}
\E[Y | U, a, X] &= \int h(w, a, X) dF(w | U, X), a = 0, 1 \label{eq: h-int}\\
\frac{1}{\PP[A = a | U, X]} &= \int q(z, a, X) dF(z | U, a, X), a = 0, 1 \label{eq: q-int}
\end{align}

\paragraph{Coefficients of $h$:}

We have that $\E[Y| U, A, X] = \gamma_0 + \gamma_a A + \gamma_x X + \gamma_u U + \gamma_{au} AU$, so \eqref{eq: h-int} implies
\small
\begin{gather*}
    \gamma_0 + \gamma_a A + \gamma_x X + \gamma_u U + \gamma_{au} AU = b_0 + b_a A + b_x X + b_{ax} AX + \int (b_w + b_{aw} A) w \cdot dF(w | U, X) \iff \\
    \gamma_0 + \gamma_a A + \gamma_x X + \gamma_u U + \gamma_{au} AU = b_0 + b_a A + b_x X + b_{ax} AX + (b_w + b_{aw} A) \E[W | U, X]
\end{gather*}
\normalsize
Since $W | U, X \sim \mathcal{N}(\mu_0 + \mu_x X + \mu_{u} U, 1)$, we get
\small
\begin{gather*}
    \gamma_0 + \gamma_a A + \gamma_x X + \gamma_u U + \gamma_{au} AU = b_0 + b_a A + b_x X + b_{ax} AX + (b_w + b_{aw} A) (\mu_0 + \mu_x X + \mu_{u} U)
\end{gather*}
\normalsize
Assigning values $A = 0, 1$, we obtain the following system
\begin{align}
    0 &= \gamma_0 - b_0 - b_w \mu_0 + (\gamma_x - b_x - \mu_x b_w) X + (\gamma_u - b_w \mu_{u}) U \label{eq: h-coef-1} \\
    \begin{split} \label{eq: h-coef-2}
        0 &= (\gamma_0 + \gamma_a) - (b_0 + b_a) - (b_w + b_{aw}) \mu_0  + \\
    &+ (\gamma_x - b_x - b_{ax} - \mu_x (b_w + b_{aw})) X + (\gamma_u + \gamma_{au} - (b_w + b_{aw}) \mu_{u}) U
    \end{split}
\end{align}

Multiplying \eqref{eq: h-coef-1} by $U$ and $X$ and taking the expectation in each resulting equations yields
\begin{align*}
    0 &= \rho (\gamma_x - b_x - \mu_x b_w) + (\gamma_u - b_w \mu_u) \\
    0 &= (\gamma_x - b_x - \mu_x b_w)  + \rho (\gamma_u - b_w \mu_u)
\end{align*}
Since $\rho \in (-1, 1)$, we obtain $\gamma_x - b_x - \mu_x b_w = \gamma_u - b_w \mu_u = 0$. From \eqref{eq: h-coef-1}, this additionally implies $\gamma_0 - b_0 - b_w \mu_0 = 0$.

Similarly, from \eqref{eq: h-coef-2} we obtain $\gamma_a - b_a - \mu_0 b_{aw} = - b_{ax} - \mu_x b_{aw} = \gamma_{au} - b_{aw} \mu_u = 0$. Solving for the coefficients of $h$, we obtain the unique solution
\begin{equation*}
    (b_0, b_a, b_x, b_w, b_{ax}, b_{aw}) = \left(\gamma_0 - \frac{\mu_0 \gamma_u}{\mu_u}, \gamma_a - \frac{\mu_0 \gamma_{au}}{\mu_u}, \gamma_x - \frac{\mu_x \gamma_u}{\mu_u}, \frac{\gamma_u}{\mu_u}, -\frac{\mu_x \gamma_{au}}{\mu_u}, \frac{\gamma_{au}}{\mu_u} \right)
\end{equation*}

\paragraph{Coefficients of $q$:} We have that $\PP[A | X, U] = \frac{1}{1 + \exp\left\{(-1)^{A} (\alpha_0 + \alpha_x X + \alpha_u U)\right\}}$, such that \eqref{eq: q-int} implies
\small
\begin{gather*}
    1 + \exp\{(-1)^A (\alpha_0 + \alpha_x X + \alpha_u U)\} = \\
    = 1 + \exp\{(-1)^{1-A} (t_0 + t_a A + t_x X) \} \int \exp\{(-1)^{1-A} t_z Z \} dF(z | U, A, X)
\end{gather*}
\normalsize
Since $Z | U, A, X \sim \mathcal{N} (\theta_0+\theta_a A + \theta_u U + \theta_x X, 1)$, we get
\small
\begin{align*}
    \quad 1 &+ \exp\{(-1)^A (\alpha_0 + \alpha_x X + \alpha_u U)\} = \\
    &= 1 + \exp\{(-1)^{1-A} (t_0 + t_a A + t_x X) \} \int \exp\{(-1)^{1-A} t_z Z \} dF(z | U, A, X)  \\
    &= 1 + \exp\{(-1)^{1-A} (t_0 + t_a A + t_x X) \} \int \frac{1}{\sqrt{2\pi}} \exp\{(-1)^{1-A} t_z Z + 0.5 (Z - \theta_0 - \theta_a A - \theta_u U - \theta_x X)^2\} \\
    &= 1 + \exp \left\{(-1)^{1-A} (t_0 + t_a A + t_x X) + (-1)^{1-A} t_z (\theta_0 + \theta_a A + \theta_u U + \theta_x X) + \frac{t_z^2}{2} \right\}
\end{align*}
\normalsize
for each $A = 0,1$. This is equivalent to
\small
\begin{gather*}
    (-1)^A (\alpha_0 + \alpha_x X + \alpha_u U) = (-1)^{1-A} (t_0 + t_a A + t_x X) + (-1)^{1-A} t_z (\theta_0 + \theta_a A + \theta_u U + \theta_x X) + 0.5 t_z^2
\end{gather*}
\normalsize

Assigning values $A = 0, 1$, we obtain the system
\begin{align}
    0 &= \alpha_0 + t_0 + \theta_0 t_z - 0.5 t_z^2 + (\alpha_x + t_x + \theta_x t_z) X + (\alpha_u + \theta_u t_z) U \\ 
    0 &= \alpha_0 + (t_0 + t_a) + (\theta_0 + \theta_a) t_z + 0.5 t_z^2 + (\alpha_x + t_x + \theta_x t_z) X + (\alpha_u + \theta_u t_z) U
\end{align}
As in the outcome bridge function case, it follows that the coefficients of 1 (the constant term), $X$, and $U$ must be identically 0. We then obtain $\alpha_0 + t_0 + \theta_0 t_z - 0.5 t_z^2 = t_a + \theta_a t_z + t_z^2 = \alpha_x + t_x + \theta_x t_z = \alpha_u + \theta_u t_z = 0$,
which yields the unique solution
\begin{align*}
    (t_0, t_a, t_x, t_z) = \left(-\alpha_0 + \frac{\theta_0}{\theta_u} \alpha_u + \frac{0.5}{\theta_u^2} \alpha_u^2, -\frac{1}{\theta_u^2} \alpha_u^2 + \frac{\theta_a}{\theta_u} \alpha_u, \frac{\theta_x}{\theta_u} \alpha_u - \alpha_x, -\frac{\alpha_u}{\theta_u} \right) 
\end{align*}

\section{Proving violations of Completeness Assumption \ref{assumption: Completeness}(a)}

\subsection{Violation of Assumption \ref{assumption: Completeness}(a) under setup \eqref{eq: U2-AYZW0}} \label{sec: U2-ZW Completeness pf}

We will prove that completeness assumption \ref{assumption: Completeness}(a) is violated under the DGP \eqref{eq: U2-AYZW0} with $\boldsymbol{\theta}_u = \begin{pmatrix} \theta_{u1} & \theta_{u2}\end{pmatrix}^T$, $\theta_{u1} \ne 0$. We note that case $\theta_{u2} \ne 0$ can be treated symmetrically, by appropriately exchanging $u1$ and $u2$ in the following computations. 

For any values $u,z,a,x$, we have that
\small
\begin{align*}
\PP[U=u|Z=z,A=a,X=x] &= \frac{\PP[U=u,Z=z|A=a,X=x]}{\PP[Z=z|A=a,X=x]} \\
&= \frac{\PP[Z=z|U=u,A=a,X=x] \PP[U=u|A=a,X=x]}{\PP[Z=z|A=a,X=x]} \\
&= \frac{\PP[\epsilon_1 = z - \theta_0 - \theta_a a - \theta_x x - \theta_{u}^T u] \frac{\PP[A=a|U=u,X=x] \PP[U=u|X=x]}{\PP[A=a|X=x]}}{\PP[Z=z|A=a,X=x]} 
\end{align*}
\normalsize

Using 
\begin{align*}
    \PP[U=u|X=x] = \frac{\exp\left(\frac{\left(\rho_2 u_1 - \rho_1 u_2 \right)^2 - \left(u_2 - \rho_2 x \right)^2 - \left(u_1 - \rho_1 x \right)^2}{2 \left(1 - \rho_1^2 - \rho_2^2 \right)} \right)}{2\pi \sqrt{1-\rho_1^2-\rho_2^2}}
\end{align*}
\normalsize we get
\small
\begin{gather*}
    \PP[U=u|Z=z,A=a,X=x] = \frac{\frac{1}{\sqrt{2\pi}} \exp \left\{-0.5 \left(Z- \theta_0 - \theta_a a - \theta_x x - \theta_{u}^T u \right)^2 \right\}}{\PP[Z=z,A=a|X=x]  \left(1 + \exp \left\{\left(-1\right)^a\left(\alpha_0 + \alpha_x x + \alpha_u^T u \right) \right\}\right)} \cdot \\
    \cdot \frac{1}{2\pi \sqrt{1-\rho_1^2-\rho_2^2}} \exp\left\{\frac{\left(\rho_2 u_1 - \rho_1 u_2 \right)^2 - \left(u_2 - \rho_2 x \right)^2 - \left(u_1 - \rho_1 x \right)^2}{2 \left(1 - \rho_1^2 - \rho_2^2 \right)} \right\} = \\
    = \frac{\exp\left\{\frac{1}{2}\frac{\left(\rho_2 u_1 - \rho_1 u_2 \right)^2 - \left(u_2 - \rho_2 x \right)^2 - \left(u_1 - \rho_1 x \right)^2 }{\left(1 - \rho_1^2 - \rho_2^2 \right)} - \frac{1}{2}\left(Z- \theta_0 - \theta_a a - \theta_x x - \theta_{u}^T u \right)^2 \right\}}{\left(2\pi \right)^{3/2} \sqrt{1-\rho_1^2 - \rho_2^2} \PP[Z=z,A=a|X=x]  \left(1 + \exp \left\{\left(-1\right)^a\left(\alpha_0 + \alpha_x x + \alpha_u^T u \right) \right\}\right)}
\end{gather*}
\normalsize

Let us consider
\begin{equation}
\begin{split}
    g(U) &= u_2 \left(u_2^2 - 3 - \alpha_{u2}^2 - \frac{\alpha_{u1}^2 \theta_{u2}^2}{\theta_{u1}^2} + \frac{2 \alpha_{u1} \alpha_{u2} \theta_{u2}}{\theta_{u1}}\right) \exp\left(-\frac{u_2^2}{2} \right) \cdot \\
    &\cdot \exp\left(-\frac{\left(\rho_2 u_1 - \rho_1 u_2 \right)^2 - \left(u_2 - \rho_2 x \right)^2 - \left(u_1 - \rho_1 x \right)^2}{2 \left(1 - \rho_1^2 - \rho_2^2 \right)} \right) \cdot \\
    &\cdot \left(2 + \exp \left(-\alpha_0 -\alpha_x x - \alpha_u^T u \right) + \exp \left(\alpha_0 +\alpha_x x + \alpha_u^T u \right)\right)
    \end{split}
\end{equation}

We will prove that $\E[g(U)|Z=z,A=a,X=x] = 0$ for any values $z,a,x$. We have 
\small
\begin{gather*}
\E[g(U)|Z=z,A=a,X=x] = \int_{(-\infty,\infty)^2} g(u) \PP[U = u|Z=z,A=a,X=x] du_1 du_2 \\
= \frac{1}{\left(2\pi \right)^{3/2} \sqrt{1-\rho_1^2 - \rho_2^2} \PP[Z=z,A=a|X=x]} \int_{(-\infty,\infty)^2} u_2 \left(u_2^2 - 3 - \alpha_{u2}^2 - \frac{\alpha_{u1}^2 \theta_{u2}^2}{\theta_{u1}^2} + \frac{2 \alpha_{u1} \alpha_{u2} \theta_{u2}}{\theta_{u1}}\right) \cdot \\
\cdot \exp \left\{-\frac{1}{2}\left(Z- \theta_0 - \theta_a a - \theta_x x - \theta_{u}^T u \right)^2 - \frac{u_2^2}{2} \right\} \left(1 + \exp\left\{\left(-1\right)^{1-a}\left(\alpha_0 + \alpha_x x + \alpha_u^T u \right) \right\} \right) du_1 du_2 
\end{gather*} 
\normalsize
Let
\small
\begin{align*}
    T_1 &= \int_{(-\infty,\infty)^2} u_2 \exp \left\{-\frac{1}{2}\left(Z- \theta_0 - \theta_a a - \theta_x x - \theta_{u}^T u \right)^2 - \frac{1}{2} u_2^2 \right\} \cdot \\
    &\cdot\left(1 + \exp\left\{\left(-1\right)^{1-a}\left(\alpha_0 + \alpha_x x + \alpha_u^T u \right) \right\} \right) du_1 du_2 \\
    T_2 &= \int_{(-\infty,\infty)^2} u_2^3 \exp \left\{-\frac{1}{2}\left(Z- \theta_0 - \theta_a a - \theta_x x - \theta_{u}^T u \right)^2 - \frac{1}{2} u_2^2 \right\} \cdot \\
    &\cdot\left(1 + \exp\left\{\left(-1\right)^{1-a}\left(\alpha_0 + \alpha_x x + \alpha_u^T u \right) \right\} \right) du_1 du_2 
\end{align*}
\normalsize
We have that
\smaller
\begin{gather*}
    \int_{-\infty}^\infty \exp \left\{-\frac{1}{2}\left(Z- \theta_0 - \theta_a a - \theta_x x - \theta_{u}^T u \right)^2 \right\} \left(1 + \exp\left\{\left(-1\right)^{1-a}\left(\alpha_0 + \alpha_x x + \alpha_u^T u \right) \right\} \right) du_1 = \\
    = \exp\left\{-\frac{1}{2} \left(Z - \theta_0 - \theta_a a - \theta_x x - \theta_{u2} u_2 \right)^2 \right\} \left(\int_{-\infty}^\infty  \exp \left\{\theta_{u1} \left(Z - \theta_0 - \theta_a a - \theta_x x - \theta_{u2} u_2 \right) u_1 - \frac{1}{2} \theta_{u1}^2 u_1^2 \right\} du_1 \right.+ \\ \left. + \exp\left\{(-1)^{1-a} \left(\alpha_0 + \alpha_x x + \alpha_{u2} u_2 \right) \right\} \cdot \right. \\
    \left. \cdot \int_{-\infty}^\infty \exp \left\{\theta_{u1} \left(Z - \theta_0 - \theta_a a - \theta_x x - \theta_{u2} u_2 - \frac{(-1)^{a} \alpha_{u1}}{\theta_{u1}}\right) u_1 - \frac{1}{2} \theta_{u1}^2 u_1^2 \right\} du_1 \right) = \\
    = \exp\left\{-\frac{1}{2} \left(Z - \theta_0 - \theta_a a - \theta_x x - \theta_{u2} u_2 \right)^2 \right\} \left( \frac{\sqrt{2\pi}}{|\theta_{u1}|} \exp \left\{\frac{\theta_{u1}^2 \left(Z - \theta_0 - \theta_a a - \theta_x x - \theta_{u2} u_2 \right)^2}{2 \theta_{u1}^2} \right\} \right. + \\
    + \exp\left\{(-1)^{1-a} \left(\alpha_0 + \alpha_x x + \alpha_{u2} u_2 \right) \right\} \cdot \left. \frac{\sqrt{2\pi}}{|\theta_{u1}|} \exp \left\{\frac{\theta_{u1}^2 \left(Z - \theta_0 - \theta_a a - \theta_x x - \theta_{u2} u_2 - \frac{(-1)^a \alpha_{u1}}{\theta_{u1}} \right)^2}{2 \theta_{u1}^2} \right\} \right) = \\
    = \frac{\sqrt{2\pi}}{|\theta_{u1}|} \left(1 + \exp\left\{(-1)^{1-a} \left(\alpha_0 + \alpha_x x + \alpha_{u2} u_2 + \frac{ \alpha_{u1}}{\theta_{u1}} \left(Z - \theta_0 - \theta_a a - \theta_x x - \theta_{u2} u_2 \right) \right) + \frac{\alpha_{u1}^2}{2\theta_{u1}^2} \right\}  \right) = \\
    = \frac{\sqrt{2\pi}}{|\theta_{u1}|} \left(1 + \exp\left\{(-1)^{1-a} \left(\alpha_0 + \alpha_x x + \frac{ \alpha_{u1}}{\theta_{u1}} \left(Z - \theta_0 - \theta_a a - \theta_x x \right) \right) + \frac{\alpha_{u1}^2}{2\theta_{u1}^2} + (-1)^a \left(\frac{\alpha_{u1} \theta_{u2}}{\theta_{u1}} - \alpha_{u2} \right) u_2 \right\}  \right)
\end{gather*}
\normalsize
which implies
\small
\begin{gather*}
    T_1 = \frac{\sqrt{2\pi}}{|\theta_{u1}|} \exp\left\{(-1)^{1-a} \left(\alpha_0 + \alpha_x x + \frac{ \alpha_{u1}}{\theta_{u1}} \left(Z - \theta_0 - \theta_a a - \theta_x x \right) \right) + \frac{\alpha_{u1}^2}{2\theta_{u1}^2} \right\} \cdot \\ \cdot \int_{-\infty}^\infty u_2 \exp \left\{(-1)^a \left(\frac{\alpha_{u1} \theta_{u2}}{\theta_{u1}} - \alpha_{u2} \right) u_2 - \frac{1}{2} u_2^2 \right\} du_2 = \\
    = \frac{\sqrt{2\pi}}{|\theta_{u1}|} \exp\left\{(-1)^{1-a} \left(\alpha_0 + \alpha_x x + \frac{ \alpha_{u1}}{\theta_{u1}} \left(Z - \theta_0 - \theta_a a - \theta_x x \right) \right) + \frac{\alpha_{u1}^2}{2\theta_{u1}^2} \right\} \cdot \\
    \cdot \sqrt{2\pi} (-1)^a \left(\frac{\alpha_{u1} \theta_{u2}}{\theta_{u1}} - \alpha_{u2} \right) \exp \left\{\frac{1}{2}\left(\frac{\alpha_{u1} \theta_{u2}}{\theta_{u1}} - \alpha_{u2} \right)^2 \right\} = \\
    = \frac{2\pi (-1)^a}{|\theta_{u1}|} \exp\left\{(-1)^{1-a} \left(\alpha_0 + \alpha_x x + \frac{ \alpha_{u1}}{\theta_{u1}} \left(Z - \theta_0 - \theta_a a - \theta_x x \right) \right) \right. + \\
    + \left. \frac{\alpha_{u1}^2 (1+\theta_{u2}^2)}{2\theta_{u1}^2} - \frac{\theta_{u1}\alpha_{u1} \alpha_{u2}}{\theta_{u1}} - \frac{1}{2} \alpha_{u2}^2 \right\} \cdot \left(\frac{\alpha_{u1} \theta_{u2}}{\theta_{u1}} - \alpha_{u2} \right)
\end{gather*}
\normalsize
and 
\small
\begin{gather*}
    T_2 = \frac{\sqrt{2\pi}}{|\theta_{u1}|} \exp\left\{(-1)^{1-a} \left(\alpha_0 + \alpha_x x + \frac{ \alpha_{u1}}{\theta_{u1}} \left(Z - \theta_0 - \theta_a a - \theta_x x \right) \right) + \frac{\alpha_{u1}^2}{2\theta_{u1}^2} \right\} \cdot \\ \cdot \int_{-\infty}^\infty u_2^3 \exp \left\{(-1)^a \left(\frac{\alpha_{u1} \theta_{u2}}{\theta_{u1}} - \alpha_{u2} \right) u_2 - \frac{1}{2} u_2^2 \right\} du_2 = \\
    = \frac{\sqrt{2\pi}}{|\theta_{u1}|} \exp\left\{(-1)^{1-a} \left(\alpha_0 + \alpha_x x + \frac{ \alpha_{u1}}{\theta_{u1}} \left(Z - \theta_0 - \theta_a a - \theta_x x \right) \right) + \frac{\alpha_{u1}^2}{2\theta_{u1}^2} \right\} \cdot \\
    \cdot \sqrt{2\pi} (-1)^a \left(\frac{\alpha_{u1} \theta_{u2}}{\theta_{u1}} - \alpha_{u2} \right) \left[3 + \left(\frac{\alpha_{u1} \theta_{u2}}{\theta_{u1}} - \alpha_{u2} \right)^2 \right] \exp \left\{\frac{1}{2}\left(\frac{\alpha_{u1} \theta_{u2}}{\theta_{u1}} - \alpha_{u2} \right)^2 \right\} = \\
    = \frac{2\pi (-1)^a}{|\theta_{u1}|} \exp\left\{(-1)^{1-a} \left(\alpha_0 + \alpha_x x + \frac{ \alpha_{u1}}{\theta_{u1}} \left(Z - \theta_0 - \theta_a a - \theta_x x \right) \right) \right. + \\
    + \left. \frac{\alpha_{u1}^2 (1+\theta_{u2}^2)}{2\theta_{u1}^2} - \frac{\theta_{u1}\alpha_{u1} \alpha_{u2}}{\theta_{u1}} - \frac{1}{2} \alpha_{u2}^2 \right\} \cdot \left(\frac{\alpha_{u1} \theta_{u2}}{\theta_{u1}} - \alpha_{u2} \right) \left(3 + \alpha_{u2}^2 + \frac{\alpha_{u1}^2 \theta_{u2}^2}{\theta_{u1}^2} - \frac{2 \alpha_{u1} \alpha_{u2} \theta_{u2}}{\theta_{u1}} \right)
\end{gather*}
\normalsize using the fact that $\int_{-\infty}^\infty u_2 \exp\left\{-\frac{1}{2} u_2^2 \right\} du_2 = 0$ and $\int_{-\infty}^\infty u_2^3 \exp\left\{-\frac{1}{2} u_2^2 \right\} du_2 = 0$ (as integrals of odd functions). We then obtain
\begin{gather*}
    \E[g(U)|Z=z,A=a,X=x] = \frac{1}{\left(2\pi \right)^{3/2} \sqrt{1-\rho_1^2 - \rho_2^2} \PP[Z=z,A=a|X=x]} \cdot \\
    \cdot \left(T_2 - \left(3 + \alpha_{u2}^2 + \frac{\alpha_{u1}^2 \theta_{u2}^2}{\theta_{u1}^2} - \frac{2 \alpha_{u1} \alpha_{u2} \theta_{u2}}{\theta_{u1}} \right) T_1 \right) = 0
\end{gather*}
for any $z,a,x$. However, we clearly do not have $g(U) \equiv 0$ a.s., so completeness assumption \ref{assumption: Completeness}(a) does not hold.  

\section{Bias computations}

\subsection{Computing the (asymptotic) bias obtained through Method of Moments estimator under setup \eqref{eq: U2-AY}}\label{sec: bias UW-AY pf}

We will compute the asymptotic bias obtained from the method of moments solver using bridge function $h(W, A,0;b) = b_0 + b_a A + b_w W + b_{aw} AW$ and vector function $Q(A, Z, 0) = (1, A, Z, AZ)^T$. 

We define the moment restrictions $H(D_i;\theta) = \begin{Bmatrix} \{Y_i - h(W_i, A_i, 0; b)\} \times Q(A_i, Z_i, 0) \\ \Delta - (h(W_i, 1, 0; b) - h(W_i, 0, 0; b)) \end{Bmatrix}$ and let $m(\theta) = \E[H(D; \theta)] = \lim_{n \rightarrow \infty} \frac{1}{n} \sum_{i=1}^n h(D_i; \theta)$. The estimate of $\theta = (b, \Delta)$ is given by
\begin{equation*}
    \hat \theta = \argmin_{\theta} m^T(\theta) m(\theta)
\end{equation*}

\small
Using $\E[U_1] = \E[U_2] = 0$, $\E[U_1^2] = \E[U_2^2] = 1$, and $\E[U_1 U_2] = 0$, we express the coordinates of $\E[h(D; \theta)] = (m_1, m_2, m_3, m_4, m_5)$ as follows:
\begin{equation}
\begin{split}
m_1 &= -b_0 - \E[A] b_a - \mu_0 b_w - (\E[A] \mu_0 + \E[AU_1] \mu_{u1}) b_{aw} + \\ &+ \gamma_0 + \E[A] \gamma_a + \E[AU_1] \gamma_{au1} 
\end{split}
\end{equation}
\begin{equation}
\begin{split}
m_2 &=  -\E[A] b_0 - \E[A] b_a - (\E[A] \mu_0 + \E[AU_1] \mu_{u1}) b_w - (\E[A] \mu_0 + \E[AU_1] \mu_{u1}) b_{aw} + \\
&+ \left(\E[A] (\gamma_0 + \gamma_a) + \E[AU_1] (\gamma_{u1} + \gamma_{au1}) \right)
\end{split}
\end{equation}

\begin{equation}
\begin{split}
m_3 &= -(\theta_0 + \E[A] \theta_a) b_0 - \left(\E[A] (\theta_0 + \theta_a) + \E[AU_1] \theta_{u1} \right) b_a - \\
&- (\mu_0 \theta_0 + \mu_{u1} \theta_{u1} + \E[A] \mu_0 \theta_a + \E[AU_1] \mu_{u1} \theta_a) b_w - \\
&- \left(\E[A] \mu_0 (\theta_0 + \theta_a) + \E[AU_1] \left(\mu_0 \theta_{u1} + \mu_{u1} (\theta_0 + \theta_a)\right) + \E[AU_1^2] \mu_{u1} \theta_{u1} \right) b_{aw} + \\
&+ \gamma_0 \theta_0 + \gamma_{u1} \theta_{u1} + \E[A] \left(\gamma_0 \theta_a + \gamma_a (\theta_0 + \theta_a) \right) + \E[AU_1] \left(\gamma_a \theta_{u1} + \gamma_{u1} \theta_a + \gamma_{au1} (\theta_0 + \theta_a) \right) + \\
&+ \E[AU_1^2] \gamma_{au1} \theta_{u1} + \E[AU_2] \gamma_{u2} \theta_a
\end{split} 
\end{equation}

\begin{equation}
\begin{split}
m_4 &= - \left( \E[A] (\theta_0 + \theta_a) + \E[AU] \theta_{u1} \right) b_0 - \left(\E[A] (\theta_0 + \theta_a) + \E[AU] \theta_{u1} \right) b_a - \\
&- \left(\E[A] \mu_0 (\theta_0 + \theta_a) + \E[AU] \left(\mu_0 \theta_{u1} + \mu_{u1} (\theta_0 + \theta_a) \right) + \E[AU^2] \mu_{u1} \theta_{u1} \right) b_w - \\
&- \left(\E[A] \mu_0 (\theta_0 + \theta_a)  + \E[AU] \left(\mu_0 \theta_{u1} + \mu_{u1} (\theta_0 + \theta_a) \right) + \E[AU^2] \mu_{u1} \theta_{u1} \right) b_{aw} - \\
&+ \E[A] (\gamma_0 + \gamma_a) (\theta_0 + \theta_a) + \E[AU] \left(\gamma_{u1} \theta_a + \gamma_{au1} (\theta_0 + \theta_a) \right) + \\
&+ \E[AU^2] \gamma_{au1} \theta_{u1} + \E[AU_2] \gamma_{u2} \theta_a
\end{split}
\end{equation}
\normalsize

Let 
\begin{align}
    R_1 &= \frac{\E[AU_1] \E[AU_2]}{(1-\E[A])(1-\E[AU_1^2]) - \E[AU_1]^2} \\
    R_2 &= \frac{1 - \E[A] - \E[AU_1^2]}{\E[A] \E[AU_1^2] - \E[AU_1]^2}
\end{align}

We obtain the estimated bridge function parameters
\begin{equation}
    \begin{split}
        \hat b_0 &= \gamma_0 + \frac{\mu_0}{\mu_{u1}} \gamma_{u1} + \left(\frac{\mu_0}{\mu_{u1}} - \frac{1 - \E[AU_1^2]}{\E[AU_1]} \right)R_1 \cdot \gamma_{u2} \\
        \hat b_a &= \gamma_a - \frac{\mu_0}{\mu_{u1}} \gamma_{au1} - \left(\frac{\mu_0}{\mu_{u1}} + \frac{\E[AU_1^2]}{\E[AU_1]} + \frac{\frac{\E[A] \E[AU_1^2]}{\E[AU_1]} - \E[AU_1]}{1 - \E[A] - \E[AU_1^2]} \right) R_1 R_2 \cdot \gamma_{u2}\\
        \hat b_w &= \frac{1}{\mu_{u1}} + \frac{R_1}{\mu_{u1}} \gamma_{u2}\\
        \hat b_{aw} &= \frac{1}{\mu_{u1}} \gamma_{au1} + \frac{R_1 R_2}{\mu_{u1}} \gamma_{u2}
    \end{split}
\end{equation}

The estimated effect resulting from $\hat h(W,A,0;b)$ is then
\small
\begin{align*}
    \hat \Delta &= \hat b_a + \hat b_{aw} \E[W] = \hat b_a + \hat b_{aw} \mu_0 =\\
    &= \gamma_a - \frac{\mu_0}{\mu_{u1}} \gamma_{au1} - \left(\frac{\mu_0}{\mu_{u1}} + \frac{\E[AU_1^2]}{\E[AU_1]} + \frac{\frac{\E[A] \E[AU_1^2]}{\E[AU_1]} - \E[AU_1]}{1 - \E[A] - \E[AU_1^2]} \right) R_1 R_2 \cdot \gamma_{u2} + \frac{\mu_0}{\mu_{u1}} \gamma_{au1} + \frac{\mu_0}{\mu_{u1}} R_1 R_2 \gamma_{u2} \\
    &= \gamma_a - \left(\frac{\E[AU_1^2]}{\E[AU_1]} + \frac{\frac{\E[A] \E[AU_1^2]}{\E[AU_1]} - \E[AU_1]}{1 - \E[A] - \E[AU_1^2]} \right) R_1 R_2 \cdot \gamma_{u2} 
\end{align*}
\normalsize
which yields a bias equal to
\begin{equation}
    \delta = - \frac{\E[AU_1^2]}{\E[AU_1]}\left(1 + \frac{\E[A]\E[AU_1^2] - \E[AU_1]^2}{\E[AU_1^2] \left(1 - \E[A] - \E[AU_1^2]\right)} \right) R_1 R_2 \cdot \gamma_{u2} 
\end{equation}

We note that the expectations
\begin{equation}
    \begin{split}
        \E[A] &= \E[\E[A|U_1, U_2]] = \E\left[\PP[A=1|U_1, U_2]\right] = \\
        &= \E\left[\frac{1}{1+\exp\{-\alpha_0 - \alpha_{u1} U_1 - \alpha_{u2} U_2\}} \right] = \int_{-\infty}^\infty \int_{-\infty}^\infty \frac{\frac{1}{2\pi} \exp\left\{-\frac{u^2 + v^2}{2}\right\} du dv}{1 + \exp\{-\alpha_0 - \alpha_{u1} u - \alpha_{u2} v\}} 
    \end{split}
\end{equation}
\begin{equation}
    \begin{split}
        \E[AU_1] &= \E[\E[AU_1|U_1, U_2]] = \E[U_1\E[A|U_1, U_2]] = \\
        &= \E\left[\frac{U_1}{1+\exp\{-\alpha_0 - \alpha_{u1} U_1 - \alpha_{u2} U_2\}} \right] = \int_{-\infty}^\infty \int_{-\infty}^\infty \frac{\frac{1}{2\pi} u \exp\left\{-\frac{u^2 + v^2}{2}\right\} du dv}{1 + \exp\{-\alpha_0 - \alpha_{u1} u - \alpha_{u2} v\}} 
    \end{split}
\end{equation}
\begin{equation}
    \begin{split}
        \E[AU_2] &= \E[\E[AU_2|U_1, U_2]] = \E[U_2\E[A|U_1, U_2]] = \\
        &= \E\left[\frac{U_2}{1+\exp\{-\alpha_0 - \alpha_{u1} U_1 - \alpha_{u2} U_2\}} \right] = \int_{-\infty}^\infty \int_{-\infty}^\infty \frac{\frac{1}{2\pi} v \exp\left\{-\frac{u^2 + v^2}{2}\right\} du dv}{1 + \exp\{-\alpha_0 - \alpha_{u1} u - \alpha_{u2} v\}} 
    \end{split}
\end{equation}
\begin{equation}
    \begin{split}
        \E\left[AU_1^2 \right] &= \E\left[\E \left[AU_1^2|U_1, U_2 \right] \right] = \E \left[U_1^2\E \left[A|U_1, U_2 \right] \right] = \\
        &= \E\left[\frac{U_1^2}{1+\exp\{-\alpha_0 - \alpha_{u1} U_1 - \alpha_{u2} U_2\}} \right] = \int_{-\infty}^\infty \frac{\frac{1}{2\pi} u^2 \exp\left\{-\frac{u^2 + v^2}{2}\right\} du dv}{1 + \exp\{-\alpha_0 - \alpha_{u1} u - \alpha_{u2} v\}} 
    \end{split}
\end{equation}
cannot be computed in closed form but can be obtained numerically using a software like Mathematica or Maple once we provide the values of $\alpha_0$ and $\alpha_u$.

\subsection{Computing the (asymptotic) bias obtained through Method of Moments estimator under setup \eqref{eq: U2-ZW}}\label{sec: bias UW-ZW pf}

We will compute the asymptotic bias obtained from the method of moments solver using bridge function $h(W, A,0;b) = b_0 + b_a A + b_w W + b_{aw} AW$ and vector function $Q(A, Z, 0) = (1, A, Z, AZ)^T$. 

We define the moment restrictions $H(D_i;\theta) = \begin{Bmatrix} \{Y_i - h(W_i, A_i, 0; b)\} \times Q(A_i, Z_i, 0) \\ \Delta - (h(W_i, 1, 0; b) - h(W_i, 0, 0; b)) \end{Bmatrix}$ and let $m(\theta) = \E[H(D; \theta)] = \lim_{n \rightarrow \infty} \frac{1}{n} \sum_{i=1}^n h(D_i; \theta)$. The estimate of $\theta = (b, \Delta)$ is given by
\begin{equation*}
    \hat \theta = \argmin_{\theta} m^T(\theta) m(\theta)
\end{equation*}

\small
Using $\E[U_1] = \E[U_2] = 0$, $\E[U_1^2] = \E[U_2^2] = 1$, and $\E[U_1U_2] = 0$, we express the coordinates of $\E[h(D; \theta)] = (m_1, m_2, m_3, m_4, m_5)$ as follows:
\begin{align*}
    m_1 &= -b_0 - \E[A] b_a - \mu_0 b_w - (\E[A] \mu_0 + \E[AU_1] \mu_{u1}) b_{aw} + \gamma_0 + \E[A] \gamma_a + \E[AU_1] \gamma_{au1} \\
    m_2 &=  -\E[A] b_0 - \E[A] b_a - (\E[A] \mu_0 + \E[AU] \mu_{u1}) b_w - (\E[A] \mu_0 + \E[AU_1] \mu_{u1}) b_{aw} + \\
&+ \left(\E[A] (\gamma_0 + \gamma_a) + \E[AU_1] (\gamma_{u1} + \gamma_{au1}) \right) \\
    m_3 &= -(\theta_0 + \E[A] \theta_a) b_0 - \left(\E[A] (\theta_0 + \theta_a) + \E[AU_1] \theta_{u1} \right) b_a - \\
&- (\mu_0 \theta_0 + \mu_{u1} \theta_{u1} + \mu_{u2} \theta_{u2} + \E[A] \mu_0 \theta_a + \E[AU_1] \mu_{u1} \theta_a) b_w - \\
&- \left(\E[A] \left(\mu_0 (\theta_0 + \theta_a) + \mu_{u2} \theta_{u2} \right) + \E[AU_1] \left(\mu_0 \theta_{u1} + \mu_{u1} (\theta_0 + \theta_a)\right) + \E[AU_1^2] \mu_{u1} \theta_{u1} \right) b_{aw} + \\
&+ \gamma_0 \theta_0 + \gamma_{u1} \theta_{u1} + \E[A] \left(\gamma_0 \theta_a + \gamma_a (\theta_0 + \theta_a) \right) + \E[AU_1] \left(\gamma_a \theta_{u1} + \gamma_{u1} \theta_a + \gamma_{au1} (\theta_0 + \theta_a) \right) + \\
&+ \E[AU_1^2] \gamma_{au1} \theta_{u1}
\end{align*}

\begin{equation*}
\begin{split}
m_4 &= - \left( \E[A] (\theta_0 + \theta_a) + \E[AU] \theta_{u1} \right) b_0 - \left(\E[A] (\theta_0 + \theta_a) + \E[AU_1] \theta_{u1} \right) b_a - \\
&- \left(\E[A] \left(\mu_0 (\theta_0 + \theta_a) + \mu_{u2} \theta_{u2} \right) + \E[AU_1] \left(\mu_0 \theta_{u1} + \mu_{u1} (\theta_0 + \theta_a) \right) + \E[AU_1^2] \mu_{u1} \theta_{u1} \right) b_w - \\
&- \left(\E[A] \left(\mu_0 (\theta_0 + \theta_a) + \mu_{u2} \theta_{u2} \right) + \E[AU_1] \left(\mu_0 \theta_{u1} + \mu_{u1} (\theta_0 + \theta_a) \right) + \E[AU_1^2] \mu_{u1} \theta_{u1} \right) b_{aw} - \\
&+ \E[A] (\gamma_0 + \gamma_a) (\theta_0 + \theta_a) + \E[AU_1] \left(\gamma_{u1} \theta_a + \gamma_{au1} (\theta_0 + \theta_a) \right) + \\
&+ \E[AU_1^2] \gamma_{au1} \theta_{u1}
\end{split}
\end{equation*}
\normalsize

Let 
\begin{align*}
    S_1 &= \frac{(1-\E[A])^2}{(1-\E[A])(1-\E[AU_1^2])-\E[AU_1]^2}\\
    S_2 &= \frac{\E[A]^2}{\E[A]\E[AU_1^2]-\E[AU_1]^2}
\end{align*}
We obtain the estimated bridge function parameters
\small
\begin{equation}
    \begin{split}
        \hat b_0 &= \gamma_0 - \left(\frac{\mu_0}{\mu_{u1} + S_1 \cdot \frac{\theta_{u2}}{\theta_{u1}} \mu_{u2} } + \frac{S_1}{1-\E[A]} \cdot \frac{ \frac{\theta_{u2}}{\theta_{u1}}\mu_{u2}}{\mu_{u1} + S_1 \cdot \frac{\theta_{u2}}{\theta_{u1}} \mu_{u2}} \right) \gamma_{u1} \\
        \hat b_a &= \gamma_a - \frac{\mu_0 - \frac{\E[AU_1]}{\E[A]}S_2 \cdot \frac{\theta_{u2}}{\theta_{u1}} \mu_{u2}}{\mu_{u1} + S_2 \cdot \frac{\theta_{u2}}{\theta_{u1}}\mu_{u2}} \gamma_{au1} + \left(\frac{\mu_0 + \frac{\E[AU_1]}{1-\E[A]}S_1 \cdot \frac{\theta_{u2}}{\theta_{u1}} \mu_{u2}}{\mu_{u1} + S_1 \cdot \frac{\theta_{u2}}{\theta_{u1}}\mu_{u2}} -\frac{\mu_0 - \frac{\E[AU_1]}{\E[A]}S_2 \cdot \frac{\theta_{u2}}{\theta_{u1}} \mu_{u2}}{\mu_{u1} + S_2 \cdot \frac{\theta_{u2}}{\theta_{u1}}\mu_{u2}} \right) \gamma_{u1} \\
        \hat b_w &= \frac{1}{\mu_{u1} + S_1 \cdot \frac{\theta_{u2}}{\theta_{u1}}\mu_{u2}} \gamma_{u1} \\
        \hat b_{aw} &= \frac{1}{\mu_{u1} + S_2 \cdot \frac{\theta_{u2}}{\theta_{u1}} \mu_{u2}} \gamma_{au1} - \left(\frac{1}{\mu_{u1} + S_1 \cdot \frac{\theta_{u2}}{\theta_{u1}} \mu_{u2}} - \frac{1}{\mu_{u1} + S_2 \cdot \frac{\theta_{u2}}{\theta_{u1}} \mu_{u2}} \right) \gamma_{u1}
    \end{split}
\end{equation}
\normalsize

The estimated effect resulting from $\hat h(W,A,0;b)$ is then
\smaller
\begin{align*}
    \hat \Delta &= \hat b_a + \hat b_{aw} \E[W] = \hat b_a + \hat b_{aw} \mu_0 = \\
    &= \gamma_a + \frac{\frac{\E[AU_1]}{\E[A]} S_2 \cdot \frac{\theta_{u2}}{\theta_{u1}} \mu_{u2}}{\mu_{u1} + S_2 \cdot \frac{\theta_{u2}}{\theta_{u1}} \mu_{u2}} \gamma_{au1} + \left(\frac{\frac{\E[AU_1]}{1-\E[A]}S_1 \cdot \frac{\theta_{u2}}{\theta_{u1}} \mu_{u2}}{\mu_{u1} + S_1 \cdot \frac{\theta_{u2}}{\theta_{u1}}\mu_{u2}} + \frac{\frac{\E[AU_1]}{\E[A]}S_2 \cdot \frac{\theta_{u2}}{\theta_{u1}} \mu_{u2}}{\mu_{u1} + S_2 \cdot \frac{\theta_{u2}}{\theta_{u1}}\mu_{u2}} \right) \gamma_{u1} = \\
    &= \gamma_a + \frac{\E[AU_1]}{\E[A](1-\E[A])} \frac{\theta_{u2}}{\theta_{u1}} \mu_{u2} \left[\frac{(1-\E[A]) S_2}{\mu_{u1} + S_2 \cdot \frac{\theta_{u2}}{\theta_{u1}}\mu_{u2}} \gamma_{au1} + \left(\frac{\E[A] S_1}{\mu_{u1} + S_1 \cdot \frac{\theta_{u2}}{\theta_{u1}}\mu_{u2}} + \frac{(1-\E[A]) S_2}{\mu_{u1} + S_2 \cdot \frac{\theta_{u2}}{\theta_{u1}}\mu_{u2}} \right) \gamma_{u1} \right]
\end{align*}
\normalsize
which yields a bias equal to
\small
\begin{equation}
    \delta =\frac{\E[AU_1]}{\E[A](1-\E[A])} \frac{\theta_{u2}}{\theta_{u1}} \mu_{u2} \left[\frac{(1-\E[A]) S_2}{\mu_{u1} + S_2 \cdot \frac{\theta_{u2}}{\theta_{u1}}\mu_{u2}} \gamma_{au1} + \left(\frac{\E[A] S_1}{\mu_{u1} + S_1 \cdot \frac{\theta_{u2}}{\theta_{u1}}\mu_{u2}} + \frac{(1-\E[A]) S_2}{\mu_{u1} + S_2 \cdot \frac{\theta_{u2}}{\theta_{u1}}\mu_{u2}} \right) \gamma_{u1} \right]
\end{equation}
\normalsize

In the particular case $\gamma_{au1} = 0$, we obtain a bias equal to
\begin{equation}
    \delta = \frac{\E[AU_1]}{\E[A](1-\E[A])} \frac{\theta_{u2}}{\theta_{u1}} \mu_{u2} \left(\frac{\E[A] S_1}{\mu_{u1} + S_1 \cdot \frac{\theta_{u2}}{\theta_{u1}}\mu_{u2}} + \frac{(1-\E[A]) S_2}{\mu_{u1} + S_2 \cdot \frac{\theta_{u2}}{\theta_{u1}}\mu_{u2}} \right) \gamma_{u1}
\end{equation}

Similarly to Proof \ref{sec: bias UW-AY pf}, we note that the expectations
\begin{align*}
    \E[A] &= \E[\E[A|U_1]] = \E\left[\PP[A=1|U_1]\right] = \\
        &= \E\left[\frac{1}{1+\exp\{-\alpha_0 - \alpha_{u1} U_1\}} \right] = \int_{-\infty}^\infty \frac{\frac{1}{\sqrt{2\pi}} \exp\left\{-\frac{u^2}{2}\right\}}{1 + \exp\{-\alpha_0 - \alpha_{u1} u\}} du \\
    \E[AU_1] &= \E[\E[AU_1|U_1]] = \E[U_1\E[A|U_1]] = \\
        &= \E\left[\frac{U_1}{1+\exp\{-\alpha_0 - \alpha_u U_1\}} \right] = \int_{-\infty}^\infty \frac{\frac{1}{\sqrt{2\pi}} u \exp\left\{-\frac{u^2}{2}\right\}}{1 + \exp\{-\alpha_0 - \alpha_{u1} u\}} du \\
        \E\left[AU_1^2 \right] &= \E\left[\E\left[AU_1^2|U_1 \right]\right] = \E\left[U_1^2\E\left[A|U_1\right]\right] = \\
        &= \E\left[\frac{U_1^2}{1+\exp\{-\alpha_0 - \alpha_u U_1\}} \right] = \int_{-\infty}^\infty \frac{\frac{1}{\sqrt{2\pi}} u^2 \exp\left\{-\frac{u^2}{2}\right\}}{1 + \exp\{-\alpha_0 - \alpha_{u1} u\}} du
\end{align*}
cannot be computed in closed form but can be obtained numerically using a software like Mathematica or Maple once we provide the values of $\alpha_0$ and $\alpha_u$.

\subsection{Computing the estimator bias under setup \eqref{eq: U2-AY}} \label{sec: g-bias-1}

Let $M = \begin{pmatrix} 1 & Z & W & A & AZ & AW\end{pmatrix}$. By the typical formula $\hat{\boldsymbol{b}} = \left(M^T M \right)^{-1} M^T Y$ for the OLS estimator and the following
\begin{align*}
    \E[Z] &= \theta_0 + \theta_a \E[A] \\
    \E[W] &= \mu_0 \\ 
    \E[AZ] &= (\theta_0 + \theta_a) \E[A] + \theta_{u1} \E[AU_1] \\ 
    \E[AW] &= \mu_0 \E[A] + \mu_{u1} \E[AU_1] \\
    \E[Z^2] &= \theta_0^2 + \theta_{u1}^2 + \theta_{u2}^2 + 1 + (\theta_a^2 + 2\theta_0 \theta_a) \E[A] + 2 \theta_0 \theta_{u1} \E[AU_1] \\
    \E[W^2] &= \mu_0^2 + \mu_{u1}^2 + \mu_{u2}^2 + 1 \\
    \E[AZ^2] &= \left(1 + \left(\theta_0 + \theta_a \right)^2 + \theta_{u2}^2 \right) \E[A] + 2 \left(\theta_0 + \theta_a \right) \theta_{u1} \E[AU_1] + \theta_{u1}^2 \E[AU_1^2] \\
    \E[AW^2] &= \left(1 + \mu_0^2 + \mu_{u2}^2 \right) \E[A] + 2 \mu_0 \mu_{u1} \E[AU_1] + \mu_{u1}^2 \E[AU_1^2] \\ 
    \E[AZW] &= \left(\left(\theta_0 + \theta_a \right) \mu_0 + \theta_{u2} \mu_{u2} \right) \E[A] + \left(\left(\theta_0 + \theta_a \right) \mu_{u1} + \theta_{u1} \mu_0 \right) \E[AU_1] + \theta_{u1} \mu_{u1} \E[AU_1^2] \\
    \E[Y] &= \gamma_0 + \gamma_a \E[A] + \gamma_{au1} \E[AU_1] \\
    \E[ZY] &= \theta_0 \gamma_0 + \theta_{u1} \gamma_{u1} +  \left(\left(\theta_0 + \theta_a \right) \gamma_a + \theta_a \gamma_0 \right) \E[A] + \left(\theta_{u1} \gamma_a + \theta_a \gamma_{u1} + \left(\theta_0 + \theta_a \right) \gamma_{au1} \right) \E[AU_1] + \\
    &+ \theta_{u1} \gamma_{au1} \E[AU_1^2] \\
    \E[WY] &= \mu_0 \gamma_0 + \mu_{u1} \gamma_{u1} + \mu_0 \gamma_a \E[A] + \left(\mu_0 \gamma_{au1} + \mu_{u1} \gamma_a \right) \E[AU_1] + \mu_{u1} \gamma_{au1} \E[AU_1^2] \\
    \E[AY] &= \left(\gamma_0 + \gamma_a \right) \E[A] + \left(\gamma_{u1} + \gamma_{au1} \right) \E[AU_1] \\
    \E[AZY] &= \left(\theta_0 + \theta_a \right) \left(\gamma_0 + \gamma_a \right) \E[A] + \left(\left(\theta_0 + \theta_a \right) \left(\gamma_{u1} + \gamma_{au1} \right) + \theta_{u1} \left(\gamma_0 + \gamma_a \right) \right) \E[AU_1] + \\
    &+ \theta_{u1} \left(\gamma_{u1} + \gamma_{au1} \right) \E[AU_1^2] \\
    \E[AWY] &= \mu_0 \left(\gamma_0 + \gamma_a \right) \E[A] + \left(\mu_0 \left(\gamma_{u1} + \gamma_{au1} \right) + \mu_{u1} \left(\gamma_0 + \gamma_a \right) \right) \E[AU_1] + \mu_{u1} \left(\gamma_{u1} + \gamma_{au1} \right) \E[AU_1^2]
\end{align*}
\normalsize
we obtain a linear regression estimator bias equal to
\small
\begin{equation}
    \begin{split}
        \delta_{OR} &= \frac{\left(\frac{\E[AU]}{\E[A]}-\frac{(1-\E[A]) \theta_a}{S_2} \theta_{u1} \right)\left(1+\mu_{u2}^2 \right) + \left(\frac{\E[AU]}{\E[A]}-\frac{(1-\E[A]) \theta_a}{S_2} \frac{\mu_{u1} \mu_{u2}}{\theta_{u2}} \right) \theta_{u2}^2}{\left(1+\frac{\theta_{u1}^2}{S_2} \right)\left(1+\mu_{u2}^2 \right) + \left(1+\frac{\mu_{u1}^2}{S_2} \right)\left(1+\theta_{u2}^2 \right) - \left(1 + 2 \frac{\theta_{u1} \mu_{u1} \theta_{u2} \mu_{u2}}{S_2} \right)} \gamma_{u1} \ + \\
        &+ \left[\frac{\left(1 + \theta_{u2}^2 + \mu_{u2}^2 \right) \E[AU]}{\E[A]\left(1-\E[A]\right)} \left(\left(1 + \theta_{u2}^2 + \mu_{u2}^2 \right) + \right. \right.\\
        &+ \left. \left(1 - \frac{\E[AU^2]}{\E[A] \left(1-\E[A] \right)} \right) \left(\theta_{u1}^2 \left(1+\mu_{u2}^2 \right) + \mu_{u1}^2 \left(1 + \theta_{u2}^2 \right) - 2 \theta_{u1} \mu_{u1} \theta_{u2}\mu_{u2} \right) \right) + \\
        &+ \theta_a \theta_{u1} \left(1 + \mu_{u2} - \frac{\mu_{u1} \theta_{u2} \mu_{u2}}{\theta_{u1}} \right) \left(-\left(\frac{\E[AU^2] \left(1 - \E[AU^2] \right)}{\E[A] \left(1 - \E[A] \right)} + \E[AU]^2 \left(\frac{1}{\E[A]^2 S_1} + \frac{1}{\left(1 - \E[A] \right)^2 S_2} \right) \right) \cdot \right.\\
        &\cdot \left. \left(\theta_{u1}^2 \left(1+\mu_{u2}^2 \right) + \mu_{u1}^2 \left(1 + \theta_{u2}^2 \right) - 2 \theta_{u1} \mu_{u1} \theta_{u2}\mu_{u2} \right) + \frac{\left(1 + \theta_{u2}^2 + \mu_{u2}^2 \right)}{\E[A]^2\left(1 - \E[A] \right)^2} \right. \cdot \\
        &\cdot \left. \left. \left(\E[A]^4 - \E[A]^3 (1 + 2 \E[AU^2]) + 3 \E[A]^2 (\E[AU^2] + \E[AU]^2) - \E[A] (3\E[AU]^2 + \E[AU^2]) + \E[AU^2] \right) \right) \right] \\
        &\cdot \prod_{i=1,2}\frac{1}{\left(1+\frac{\theta_{u1}^2}{S_i} \right)\left(1+\mu_{u2}^2 \right) + \left(1+\frac{\mu_{u1}^2}{S_i} \right)\left(1+\theta_{u2}^2 \right) - \left(1 + 2 \frac{\theta_{u1} \mu_{u1} \theta_{u2} \mu_{u2}}{S_i} \right)} \gamma_{au1}
    \end{split}
\end{equation}
\normalsize

In particular, for $\gamma_{au1} = 0$, we obtain a bias equal to
\begin{equation}
    \delta_{OR} = \frac{\left(\frac{\E[AU]}{\E[A]}-\frac{(1-\E[A]) \theta_a}{S_2} \theta_{u1} \right)\left(1+\mu_{u2}^2 \right) + \left(\frac{\E[AU]}{\E[A]}-\frac{(1-\E[A]) \theta_a}{S_2} \frac{\mu_{u1} \mu_{u2}}{\theta_{u2}} \right) \theta_{u2}^2}{\left(1+\frac{\theta_{u1}^2}{S_2} \right)\left(1+\mu_{u2}^2 \right) + \left(1+\frac{\mu_{u1}^2}{S_2} \right)\left(1+\theta_{u2}^2 \right) - \left(1 + 2 \frac{\theta_{u1} \mu_{u1} \theta_{u2} \mu_{u2}}{S_2} \right)} \gamma_{u1}
\end{equation}

\color{black}

\subsection{Computing the estimator bias under setup \eqref{eq: U2-ZW}}\label{sec: g-bias-2}

Let $M = \begin{pmatrix} 1 & Z & W & A & AZ & AW\end{pmatrix}$. Similarly to Appendix \ref{sec: g-bias-1}, we use the typical OLS estimator $\hat{\boldsymbol{b}} = \left(M^T M \right)^{-1} M^TY$ and the following values
\begin{align*}
    \E[Z] &= \theta_0 + \theta_a \E[A] \\
    \E[W] &= \mu_0 \\ 
    \E[AZ] &= (\theta_0 + \theta_a) \E[A] + \theta_{u1} \E[AU_1] \\ 
    \E[AW] &= \mu_0 \E[A] + \mu_{u1} \E[AU_1] \\
    \E[Z^2] &= \theta_0^2 + \theta_{u1}^2 + 1 + (\theta_a^2 + 2\theta_0 \theta_a) \E[A] + 2 \theta_0 \theta_{u1} \E[AU_1] \\
    \E[W^2] &= \mu_0^2 + \mu_{u1}^2 + 1 \\
    \E[AZ^2] &= \left(1 + \left(\theta_0 + \theta_a \right)^2 \right) \E[A] + 2 \left(\theta_0 + \theta_a \right) \theta_{u1} \E[AU_1] + \theta_{u1}^2 \E[AU_1^2] \\
    \E[AW^2] &= \left(1 + \mu_0^2 \right) \E[A] + 2 \mu_0 \mu_{u1} \E[AU_1] + \mu_{u1}^2 \E[AU_1^2] \\ 
    \E[AZW] &= \left(\theta_0 + \theta_a \right) \mu_0  \E[A] + \left(\left(\theta_0 + \theta_a \right) \mu_{u1} + \theta_{u1} \mu_0 \right) \E[AU_1] + \theta_{u1} \mu_{u1} \E[AU_1^2] \\
    \E[Y] &= \gamma_0 + \gamma_a \E[A] + \gamma_{au1} \E[AU_1] \\
    \E[ZY] &= \theta_0 \gamma_0 + \theta_{u1} \gamma_{u1} +  \left(\left(\theta_0 + \theta_a \right) \gamma_a + \theta_a \gamma_0 \right) \E[A] + \left(\theta_{u1} \gamma_a + \theta_a \gamma_{u1} + \left(\theta_0 + \theta_a \right) \gamma_{au1} \right) \E[AU_1] + \\
    &+ \theta_{u1} \gamma_{au1} \E[AU_1^2] + \theta_a \gamma_{u2} \E[AU_2]\\
    \E[WY] &= \mu_0 \gamma_0 + \mu_{u1} \gamma_{u1} + \mu_0 \gamma_a \E[A] + \left(\mu_0 \gamma_{au1} + \mu_{u1} \gamma_a \right) \E[AU_1] + \mu_{u1} \gamma_{au1} \E[AU_1^2] \\
    \E[AY] &= \left(\gamma_0 + \gamma_a \right) \E[A] + \left(\gamma_{u1} + \gamma_{au1} \right) \E[AU_1] + \gamma_{u2} \E[AU_2] \\
    \E[AZY] &= \left(\theta_0 + \theta_a \right) \left(\gamma_0 + \gamma_a \right) \E[A] + \left(\left(\theta_0 + \theta_a \right) \left(\gamma_{u1} + \gamma_{au1} \right) + \theta_{u1} \left(\gamma_0 + \gamma_a \right) \right) \E[AU_1] + \\
    &+ \theta_{u1} \left(\gamma_{u1} + \gamma_{au1} \right) \E[AU_1^2] + \left(\theta_0 + \theta_a \right) \gamma_{u2} \E[AU_2]\\
    \E[AWY] &= \mu_0 \left(\gamma_0 + \gamma_a \right) \E[A] + \left(\mu_0 \left(\gamma_{u1} + \gamma_{au1} \right) + \mu_{u1} \left(\gamma_0 + \gamma_a \right) \right) \E[AU_1] + \\
    &+ \mu_{u1} \left(\gamma_{u1} + \gamma_{au1} \right) \E[AU_1^2] + \mu_0 \gamma_{u2} \E[AU_2]
\end{align*}
to compute the regression estimator bias $\delta_{OR}$. We then obtain
\begin{equation}
    \begin{split}
        \delta_{OR} &= \left[\frac{\frac{\E[AU_1]}{\E[A]} - \frac{1-\E[A]}{S_2} \theta_a \theta_{u1}}{1 + \frac{\theta_{u1}^2 + \mu_{u1}^2}{S_2}} + \frac{\frac{\E[AU_1]}{1 - \E[A]} - \frac{\E[A]}{S_1} \theta_a \theta_{u1}}{1 + \frac{\theta_{u1}^2 + \mu_{u1}^2}{S_1}} \right] \gamma_{u1} + \\
        &+ \frac{\frac{\E[AU_1]}{\E[A]} - \frac{1-\E[A]}{S_2} \theta_a \theta_{u1}}{1 + \frac{\theta_{u1}^2 + \mu_{u1}^2}{S_2}} \gamma_{au1} + \\
        &+ \E[AU_1]\E[AU_2] \theta_a \theta_{u1} \left[\frac{1-\E[A]}{\E[A]^2 \left(1 + \frac{\theta_{u1}^2 + \mu_{u1}^2}{S_2} \right)} + \frac{\E[A]}{\left(1 -\E[A] \right)^2  \left(1 + \frac{\theta_{u1}^2 + \mu_{u1}^2}{S_1} \right)}\right] \gamma_{u2}
    \end{split}
\end{equation}
In particular, for $\gamma_{au1} = 0$, we obtain a bias equal to
\begin{equation}
\begin{split}
    \delta_{OR} &= \left[\frac{\frac{\E[AU_1]}{\E[A]} - \frac{1-\E[A]}{S_2} \theta_a \theta_{u1}}{1 + \frac{\theta_{u1}^2 + \mu_{u1}^2}{S_2}} + \frac{\frac{\E[AU_1]}{1 - \E[A]} - \frac{\E[A]}{S_1} \theta_a \theta_{u1}}{1 + \frac{\theta_{u1}^2 + \mu_{u1}^2}{S_1}} \right] \gamma_{u1} + \\
    &+ \E[AU_1]\E[AU_2] \theta_a \theta_{u1} \left[\frac{1-\E[A]}{\E[A]^2 \left(1 + \frac{\theta_{u1}^2 + \mu_{u1}^2}{S_2} \right)} + \frac{\E[A]}{\left(1 -\E[A] \right)^2  \left(1 + \frac{\theta_{u1}^2 + \mu_{u1}^2}{S_1} \right)}\right] \gamma_{u2}
\end{split}
\end{equation}

\color{black}

\subsection{Comparison of Proximal and Unadjusted Estimator Biases under Setup \eqref{eq: U2-ZW}}\label{sec: bias-comparisons-ZW}

We begin by proving that both $S_1, S_2 > 0$:

\paragraph{Proof that $S_1, S_2 > 0$:} We have that
\begin{align*}
    \left|Cov(A, AU)\right| = \left|\E \left[A^2 U \right] \right| = \left|\E[AU]\right| &\le \sqrt{Var(A) Var(AU)} = \sqrt{Var(A)} \sqrt{\E[AU^2] - \E[AU]^2} 
\end{align*}
which implies $\E[AU]^2 \le Var(A) \left(\E\left[AU^2\right] - \E[AU]^2 \right)$. It follows that
\begin{align*}
    \E[A] \E[AU^2] &\ge \E[A] \cdot \frac{\E[AU]^2 (1 + Var(A))}{Var(A)} = \E[A] \cdot \frac{\E[AU]^2 (1 + Var(A))}{\E[A] - \E[A]^2} \\ &= \frac{\E[AU]^2 (1 + Var(A))}{1-\E[A]} \\
    &\ge \E[AU]^2
\end{align*}
since $1+Var(A) \ge 1$ and $1-\E[A] \in (0,1)$. Thus, $S_2 > 0$.

Similarly, if we consider $\bar A = 1 - A$ (such that $\bar A^2 = \bar A$, $\E[\bar A] = 1-\E[A]$, $Var(\bar A) = Var(A)$, $\E[\bar A U] = -\E[AU]$, and $\E[\bar AU^2] = 1 - \E[AU^2]$), we obtain
\begin{align*}
    (1-\E[A]) (1-\E[AU^2]) = \E[\bar A] \E[\bar A U^2] \ge \E[\bar A U]^2 = \E[AU]^2
\end{align*}
Thus, $S_1 > 0$ as well. \qed

Taking the ratio of magnitudes for the two biases, we have
\begin{align*}
    \left|\frac{\delta_{POR}}{\delta_{unadj}} \right| = \left|\E[A] \cdot \frac{S_1}{\frac{\theta_{u1}\mu_{u1}}{\theta_{u2}\mu_{u2}} + S_1} + \left(1 - \E[A]\right) \cdot \frac{S_2}{\frac{\theta_{u1}\mu_{u1}}{\theta_{u2}\mu_{u2}} + S_2} \right|
\end{align*}
Let $f(r) = \E[A] \cdot \frac{S_1}{r + S_1} + \left(1 - \E[A]\right) \cdot \frac{S_2}{r + S_2}$ for $r \in \left(-\infty, -\min\{S_1, S_2\} \right)$. We note that $f(r)$ is strictly increasing in $r$, that $\lim_{r \to -\infty} f(r) = 0$, and that $f(r) = 1$ has the unique solution $r^\ast = -S_1(1-\E[A]) - S_2 \E[A] < 0$. We consider the following four cases:
\begin{itemize}
    \item [I.] If $\frac{\theta_{u1}\mu_{u1}}{\theta_{u2}\mu_{u2}} \ge 0$, then $S_1, S_2 > 0$ imply that $\frac{S_i}{\frac{\theta_{u1}\mu_{u1}}{\theta_{u2}\mu_{u2}} + S_i} \in (0,1)$ for $i = 1,2$. Since $\E[A] \in (0,1)$, it follows that $0 < \left|\frac{\hat \delta_{POR}}{\hat \delta_{unadj}} \right| < 1$.
    \item [II.] If $-\min\{S_1, S_2\} \le \frac{\theta_{u1}\mu_{u1}}{\theta_{u2}\mu_{u2}} < 0$, then $S_1, S_2 > 0$ imply that $\frac{S_i}{\frac{\theta_{u1}\mu_{u1}}{\theta_{u2}\mu_{u2}} + S_i} > 1$ for $i = 1,2$. Similarly, it follows that $\left|\frac{\hat \delta_{POR}}{\hat \delta_{unadj}} \right| > 1$. In particular, if $\frac{\theta_{u1}\mu_{u1}}{\theta_{u2}\mu_{u2}} = -\min\{S_1, S_2\}$, the proximal estimator bias can be arbitrarily large. \item [III.] If $r^\ast \le \frac{\theta_{u1}\mu_{u1}}{\theta_{u2}\mu_{u2}} < -\min\{S_1, S_2\}$, then $\left|\frac{\hat \delta_{POR}}{\hat \delta_{unadj}} \right| \ge 1$.
    \item [IV.] If $\frac{\theta_{u1}\mu_{u1}}{\theta_{u2}\mu_{u2}} < r^\ast$, then $0 \le \left|\frac{\hat \delta_{POR}}{\hat \delta_{unadj}} \right| < 1$. In particular, $\frac{\theta_{u1}\mu_{u1}}{\theta_{u2}\mu_{u2}} = -\infty$ implies that the proximal estimator is unbiased (as either $\theta_{u2} = 0$ or $\mu_{u2} = 0$).
\end{itemize}
\normalsize

\subsection{Computing the Proximal Estimator Bias under $\gamma_{au} = 0$ and $h(W,A,X) = b_0 + b_a A + b_x^T X + b_w^T W$} \label{sec: general-pi-bias}

We will compute the asymptotic bias obtained from the method of moments solver using bridge function $h(W,A,X;b) = b_0 + b_a A + b_w^T W + b_x^T X$ and vector function $q(A,Z,X) = \left(1,A,Z,X \right)$. We assume the general case of multi-dimensional $U,Z,W,X$ with $Z \in \R^m$, $W \in \R^n$, $U \in \R^p$, $X \in \R^q$. Throughout this section, we use the shorthand $\E[AU] = \left(\E[AU_1],\ldots,\E[AU_p] \right)$ and $\E[AX] = \left(\E[AX_1],\ldots,\E[AX_q] \right)$.

We define the moment restrictions $H(D_i;\theta) = \begin{Bmatrix} \{Y_i - h(W_i, A_i, X_i; b)\} \times Q(A_i, Z_i, X_i) \\ \Delta - (h(W_i, 1, X_i; b) - h(W_i, 0, X_i; b)) \end{Bmatrix}$ and let $m(\theta) = \E[H(D; \theta)] = \lim_{n \rightarrow \infty} \frac{1}{n} \sum_{i=1}^n h(D_i; \theta)$. The estimate of $\theta = (b, \Delta)$ is given by
\begin{equation*}
    \hat \theta = \argmin_{\theta} m^T(\theta) m(\theta)
\end{equation*}

Using $\E[U_i] = 0$, $\E[U_i^2] = 1, \forall i = 1,\ldots,p$, and $\E[U_i U_j] = 0, \forall i,j = 1,\ldots,p, \ i \ne j$, as well as $\E[X_j] = 0$, $\E[XX^T] = \Sigma_x$, and $\E[UX^T] = \rho$, we express the coordinates of $\E[h(D; \theta)] = (m_1, m_2, \boldsymbol{m}_3, \boldsymbol{m}_4)$ with $m_1,m_2 \in \R$, $\boldsymbol{m}_3 \in \R^m$, $\boldsymbol{m}_4 \in \R^q$ as follows:
\begin{align*}
    m_1 &= -b_0 - \E[A] b_a - \mu_0^T b_w + \gamma_0 + \E[A] \gamma_a \\
    m_2 &=  -\E[A] b_0 - \E[A] b_a - \left(\E[A] \mu_0^T + \E[AU]^T \mu_u + \E[AX]^T \mu_x \right) b_w - \E[AX]^T b_x + \\
    &+ \left(\gamma_0 + \gamma_a \right) \E[A] + \E[AU]^T \gamma_u + \E[AX]^T \gamma_x \\
    m_3 &= -(\theta_0 + \theta_a \E[A]) b_0 - \left(\E[A] (\theta_0 + \theta_a) + \theta_u^T \E[AU] + \theta_x^T \E[AX] \right) b_a - \\
&- \left( \left(\theta_0 + \E[A] \theta_a \right) \mu_0^T + \left(\theta_a \E[AU]^T + \theta_u^T + \theta_x^T \rho^T \right) \mu_u + \left(\theta_a \E[AX]^T + \theta_x^T \Sigma_x  + \theta_u^T \rho \right) \mu_x \right) b_w - \\
&- \left(\theta_a \E[AX]^T + \theta_u^T \rho + \theta_x^T \Sigma_x \right) b_x + \left(\theta_0 + \theta_a \E[A] \right) \gamma_0 + \\
&+ \left(\left(\theta_0 + \theta_a \right) \E[A] + \theta_u^T \E[AU] + \theta_x^T \E[AX] \right) \gamma_a + \\
&+ \left(\theta_a \E[AU]^T + \theta_u^T + \theta_x^T \rho^T \right)\gamma_u + \left(\theta_a \E[AX]^T + \theta_u^T \rho + \theta_x^T \Sigma_x \right) \gamma_x \\
m_4 &= -\E[AX] b_a - \left(\Sigma_x \mu_x + \rho^T \mu_u \right) b_w - \Sigma_x b_x + \E[AX] \gamma_a + \Sigma_x \gamma_x + \rho^T \gamma_u
\end{align*}

\textbf{Under assumption $m = n$ and $p > m$:} Let us define 
\begin{align*}
    \beta &= \frac{\E[AU]^T - \E[AX]^T \Sigma_x^{-1} \rho^T}{\E[A] \left(1 - \E[A] \right) - \E[AX]^T \Sigma_x^{-1} \E[AX]}\\
    B &= \left(I_p - \rho \Sigma_x^{-1} \rho^T - \frac{\left(\E[AU] - \rho \Sigma_x^{-1} \E[AX] \right) \left(\E[AU]^T - \E[AX]^T \Sigma_x^{-1} \rho^T \right)}{\E[A]\left(1-\E[A]\right) - \E[AX]^T \Sigma_x^{-1} \E[AX]} \right) \theta_u.
\end{align*}
Setting $m_1 = m_2 = m_{3i} = m_{4i} = 0$ for all $i = 1,\ldots,m$, $j = 1,\ldots,q$, we get solution
\begin{align*}
    b_0 &= \gamma_0 - \E[A] \beta \gamma_u - \left(\mu_0^T - \E[A] \beta \mu_u\right) \left(B^T \mu_u \right)^{-1} B^T \gamma_u \\
    b_a &= \gamma_a + \beta \left(I_p - \mu_u \left(B^T \mu_u \right)^{-1} B^T \right) \gamma_u \\
    b_w &= \left(B^T \mu_u \right)^{-1} B^T \gamma_u \\
    b_x &= \gamma_x + \Sigma_x^{-1} \left(\rho^T - \E[AX] \beta \right) \gamma_u - \left(\mu_x + \Sigma_x^{-1} \left(\rho^T - \E[AX] \beta \right) \mu_u \right) \left(B^T \mu_u \right)^{-1} B^T \gamma_u
\end{align*}
(assuming $B^T \mu_u$ is invertible). The estimated effect resulting from $\hat h(W,A,X;b)$ is then
\begin{align*}
    \hat \Delta = \hat b_a = \gamma_a + \frac{\E[AU]^T - \E[AX]^T \Sigma_x^{-1} \rho^T}{\E[A] \left(1 - \E[A] \right) - \E[AX]^T \Sigma_x^{-1} \E[AX]} \left(I_p - \mu_u \left(B^T \mu_u \right)^{-1} B^T \right) \gamma_u 
\end{align*}
which yields a bias equal to
\begin{equation*}
    \delta = \frac{\E[AU]^T - \E[AX]^T \Sigma_x^{-1} \rho^T}{\E[A] \left(1 - \E[A] \right) - \E[AX]^T \Sigma_x^{-1} \E[AX]} \left(I_p - \mu_u \left(B^T \mu_u \right)^{-1} B^T \right) \gamma_u 
\end{equation*}

\end{document}